\documentstyle[amssymb,amscd,12pt,a4]{amsart}
\numberwithin{equation}{section}

\newtheorem{theorem}{Theorem}[section]

\newtheorem{lemma}[theorem]{Lemma}
\newtheorem{proposition}[theorem]{Proposition}
\newtheorem{definition}{Definition}[section]
\newtheorem{remark}{Remark}[section]
\newcommand{\tsize}{\textstyle}
\newcommand{\dsize}{\displaystyle}
\newcommand{\limfunc}[1]{{\mathop{#1}\limits}}
\begin{document}
\title[Generalized q-deformed Gelfand-Dickey structures]{Generalized 
q-deformed Gelfand-Dickey structures on the group of q-pseudodifference 
operators}
\author{A.L.Pirozerski}
\author{M.A.Semenov-Tian-Shansky}
\date{}
\maketitle

\begin{center}
 Universit\'e de Bourgogne, Dijon France  \\ and 
Steklov Institute of Mathematics, St.Petersburg, Russia 
\end{center}



\begin{abstract}
We define the q-deformed Gelfand-Dickey bracket on the space of 
q-pseudodifference symbols which agrees with the Poisson Virasoro algebra of
E.Frenkel and N.Reshetikhin and its generalizations and prove its uniqueness
(in a natural class of quadratic Poisson structures). The associated 
hierarchies of nonlinear q-difference equations are also constructed.
\end{abstract}
\bigskip
\section{Introduction}

It is well known that the generalized KdV hierarchy of non-linear
differential equations admits several different realizations. The first one
is associated with the algebra of pseudodifferential operators on the line
(or on the circle). The famous construction assigns to each nonlinear
evolution equation in this hierarchy a pair $(L,A)$ of differential
operators such that the evolution equation is equivalent to the Lax equation 
\begin{equation}
\frac{dL}{dt}=\left[ A,L\right] .  \label{lax}
\end{equation}
The space of differential operators admits several remarkable Poisson
structures, and Lax equations are Hamiltonian with respect to each of them.
The simplest one is the so called first Gelfand-Dickey bracket, which is a 
{\em linear} Poisson bracket naturally arising from the identification of
the space of differential operators with the dual of the Lie algebra of
integral operators \cite{a, lm}. The next one is the celebrated second
Gelfand-Dickey bracket (or, Adler-Gelfand-Dickey bracket) \cite{a, gd}. This
bracket is {\em quadratic}, and its geometric comprehension has required
much work; it admits at least three different realization, and isomorphisms
between them usually represent deep theorems. The first one, which appears
naturally in the study of Lax equations (\ref{lax}), is based on the study
of the Lie group of integral operators (more precisely, of its central
extension \cite{KhZ}). This group comes equipped with the natural Sklyanin
bracket which endows it with the structure of a Poisson-Lie group, and the
second Gelfand-Dickey bracket is identified with the Sklyanin bracket on its
Poisson subvariety. The second realization, which is totally different, is
based on the study of the center of the universal enveloping algebra $U(%
\widehat{{\frak sl}(n)})$ of the central extension of the loop algebra of $%
{\frak sl}(n)$ at the critical value of the central charge \cite{ef}. The
third realization, finally, is naturally related to the alternative
description of the generalized KdV hierarchy which is provided by the
Drinfeld-Sokolov theory \cite{ds}. The two latter approaches provide a
natural generalization of the second Gelfand-Dickey bracket for arbitrary
semisimple Lie algebras; the corresponding Poisson algebras 
are then called classical W-algebras.

Nonlinear differential equations (\ref{lax}) admit natural difference or
q-difference analogues; their Hamiltonian treatment is more or less parallel
to the differential case, although there arise some new and unexpected
phenomena. As it happens, all three different constructions of the classical
W-algebras referred to above have their natural q-difference counterparts.
Historically, the first one to arise was based on the study of the Poisson
structure on the center of the quantized universal enveloping algebra $U_q(%
\widehat{{\frak sl}(n)})$ \cite{FrRes}. The quantization parameter $q$ is
naturally identified with the modulus of the associated q-difference
operator, $D_qf(z)=f(qz).$

The same Poisson structure also arises as a result of the Drinfeld-Sokolov
type reduction for the first order matrix q-difference equation \cite{FRS,
SemSev}. A nontrivial point in the reduction procedure is that it involves a
new elliptic classical r-matrix (its introduction is prompted by the
consistency conditions for the reduction); the modulus $\tau $ of the
underlying elliptic curve is related to $q$ via $q=\exp \pi i\tau .$

The goal of the present paper is to provide the q-difference counterpart of
the last remaining construction which is based on the study of the algebra
of q-pseudodifference symbols. We prove that for each $n\in {\Bbb N}$ there
exists a unique quadratic Poisson structure of the natural r-matrix type on
the space ${\Bbb M}_n$ of the n-th order  q-difference operators with normalized
highest term such that formal spectral invariants 
\begin{equation}
H_m(L)=\frac nm{\rm Tr\,}L^{\frac mn},\quad m\in {\Bbb N\quad },(m,n)=1,
\label{tr}
\end{equation}
of a difference operator $L=D^n+u_{n-1}D^{n-1}+\cdot \cdot \cdot +u_0$ are
in involution and generate q-difference Lax equations 
\begin{equation}
\frac{dL}{dt}=\left[ A,L\right] ,\quad A=L_{\left( +\right) }^{\frac mn};
\label{qlax}
\end{equation}
moreover, this Poisson structure coincides with the one obtained via the
q-difference Drinfeld-Sokolov reduction (or, equivalently, with the one
obtained in \cite{FrRes} via the study of the center of $U_q(\widehat{{\frak %
sl}(n)})$ at the critical level). The generalized q-deformed KdV hierarchy
which corresponds to (\ref{qlax}) was described earlier by E.Frenkel \cite
{Fr}; however, his approach to the description of the associated Poisson
structure is different: he simply uses the Poisson bracket borrowed from 
\cite{FrRes} and does not discuss its construction 
via the r-matrix formalism for the algebra of q-difference operators. Let us
also note that the lattice version of this Poisson structure has been
introduced (in a different context) by W.Oewel \cite{Oew} who also
considered the lattice analogues of the KdV and KP hierarchies. These
lattice hierarchies are also studied in \cite{Pir}.

In the second part of this paper the Poisson structure on the space of
q-difference operators is generalized to the case of q-pseudodifference
operators of arbitrary complex degree; this construction is motivated by 
\cite{KhZ}, \cite{KhLR}. We extend the algebra $\Psi {\bf D}_q$ of
q-pseudodifference symbols by adjoining to it the outer derivation $ad\ln D$
and performing the associated central extension; the extended Lie algebra of
q-integral operators gives rise to the Lie group 
\begin{equation}
\widehat{G}_{-}=\bigcup_{\alpha \in {\Bbb C}}\widehat{G}_\alpha ,\quad
\widehat{G}_\alpha =\left\{ L|\;L=D^\alpha +\sum_{i=1}^\infty
u_iD^{\alpha -i}\right\}.  \label{alpha}
\end{equation}
If $\alpha \in {\Bbb C}$ is generic, i.e., satisfies $\frac{\alpha \ln q}{%
2\pi i}\notin {\Bbb Q},$ for all elements $L\in \widehat{G}_\alpha $
there exists a logarithm and hence we may define $L^\beta $ for each $\beta
\in {\Bbb C}.$ In particular, $L^{\frac m\alpha }\in \widehat{G}_m$
for any $m\in {\Bbb N}$ and hence contains only integer powers of $D;$ let $%
L_{\left( +\right) }^{\frac m\alpha }$ be its positive part. The equation 
\begin{equation}
\frac{dL}{dt}=\left[ L_{\left( +\right) }^{\frac m\alpha },L\right] \quad
\label{f07}
\end{equation}
preserves $\widehat{G}_\alpha $ and has an infinite family of
conservation laws $H_n(L)=\frac \alpha n{\rm Tr\,}L^{\frac n\alpha},$ 
$\quad n\in {\Bbb N}.$ The flows (\ref{f07}) for different $m$ commute each with other.
We show that in a natural class of Poisson brackets on $\widehat{G}%
_\alpha $there exists a unique one with respect to which the equations (\ref
{f07}) are induced by the Hamiltonians $H_m(L)$. For integer $\alpha $ this
bracket may be restricted to ${\Bbb M}_\alpha $ ; this restriction coincides
with the bracket constructed in the first part of the present paper. A
similar class of equations has been considered in \cite{KhLR}, but Poisson
structures for them have not been proposed. We shall discuss the relation
between these two construction below (see remark \ref{rkhlr}).

\section{Nonlinear q-difference equations of the KdV type}

\subsection{Notation}

Throughout the paper we shall use the following notation. Let $\hat h$ be
the dilation operator, 
\begin{equation}
\hat hf(z)=f(qz),\quad f\in {\Bbb C}\left( \left( z^{-1}\right) \right)
,\quad q\in {\Bbb C},\quad \left| q\right| <1.  \label{f11}
\end{equation}
We denote by $\Psi {\bf D}_q$ the algebra of q-pseudodifference operators;
by definition, $\Psi {\bf D}_q$ consists of formal series of the form 
\begin{equation}
A=\sum_{i=-\infty }^{N\left( A\right) }a_i\left( z\right) D^i,\quad a_i\in 
{\Bbb C}\left( \left( z^{-1}\right) \right)  \label{f12}
\end{equation}
with the commutation relation 
\begin{equation}
D\cdot a=\left( \hat ha\right) \cdot D,\quad a\in {\Bbb C}\left( \left(
z^{-1}\right) \right) .  \label{f13}
\end{equation}
For $a\in {\Bbb C}\left( \left( z^{-1}\right) \right) $ we put 
\begin{equation}
^{h^l}a=\hat h^la,\quad \forall l\in {\Bbb C}.  \label{f14}
\end{equation}
As a linear space, $\Psi {\bf D}_q$ is a direct sum of three subalgebras, 
\begin{align}
J_{+}& =\left\{ A\in \Psi {\bf D}_q|\;A=\sum_{i=1}^{N\left( A\right)
}a_i\left( z\right) D^i,\quad a_i\in {\Bbb C}\left( \left( z^{-1}\right)
\right) \right\} ,  \label{f15} \\
J_0& ={\Bbb C}\left( \left( z^{-1}\right) \right) \subset \Psi {\bf D}_q,
\label{f16} \\
J_{-}& =\left\{ A\in \Psi {\bf D}_q|\;A=\sum_{i=1}^\infty a_i\left( z\right)
D^{-i},\quad a_i\in {\Bbb C}\left( \left( z^{-1}\right) \right) \right\} .
\label{f17}
\end{align}
Clearly, $J_0$ normalizes $J_{\pm }$ and hence $J_{\left( \pm \right)
}=J_{\pm }+J_0$ is also a subalgebra. Let $P_{\pm },P_0$ be the associated
projection operators which project $\Psi {\bf D}_q$ onto $J_{\pm },J_0,$
respectively, parallel to the complement. Put $P_{\left( \pm \right)
}=P_{\pm }+P_0.$ For $A\in \Psi {\bf D}_q$ set $A_{\pm }=$ $P_{\pm }A,\quad
A_{\left( \pm \right) }=P_{\left( \pm \right) }A,$ ${\rm \limfunc{Res}}%
A=A_0=P_0A.$ For $a\in {\Bbb C}\left( \left( z^{-1}\right) \right) ,\quad
a=\sum_ia_iz^i,$ we put 
\begin{equation}
\int a(z)dz/z=a_0;  \label{f19}
\end{equation}
clearly, this formal integral is dilation invariant, i.e., 
\begin{equation}
\int a(z)dz/z=\int a(qz)dz/z.  \label{f110}
\end{equation}
For $A\in \Psi {\bf D}_q$ we define its formal trace by 
\begin{equation}
Tr\;A=\int {\rm \limfunc{Res}}\;Adz/z;  \label{f111}
\end{equation}
it is easy to see that $Tr\;AB=Tr\;BA$ for any $A,B\in \Psi {\bf D}_q.$ We
introduce an inner product in $\Psi {\bf D}_q$ by 
\begin{equation}
\left\langle A,B\right\rangle =Tr\;AB,\quad A,B\in \Psi {\bf D}_q.
\label{f112}
\end{equation}
Clearly, this inner product is invariant and non-degenerate and the
subalgebras $J_{\pm }$ are isotropic; moreover, it sets $J_{+}$ and \ $J_{-}$
into duality, while $J_0\simeq J_0^{*}.$

\subsection{Fractional powers of q-pseudodifference operators and Lax
equations}

The fractional powers formalism which is described below is largely parallel
to the standard pseudodifferential case. Let ${\Bbb M}_{n}$ $\subset
\Psi {\bf D}_q$ be the affine subspace consisting of q-difference operators
of the form 
\begin{equation}
L=D^n+u_{n-1}D^{n-1}+\cdots +u_0,\quad u_i\in {\Bbb C}\left( \left(
z^{-1}\right) \right) .  \label{f119}
\end{equation}
We are interested in Lax equations of the form 
\begin{equation}
\frac{dL}{dt}=\left[ A,L\right] ,\quad L\in {\Bbb M}_{n}.  \label{f118}
\end{equation}
For consistency, the commutator in the r.h.s must be a polynomial in $D$ of
degree $\leq n-1.$ Let $Z_L$ be the centralizer of $L$ in $\Psi {\bf D}_q.$
Put 
\begin{equation*}
\Omega _L=\left\{ A\in \Psi {\bf D}_q|\deg \left[ A,L\right] \leq
n-1\right\} .
\end{equation*}

\begin{proposition}
\label{p11} $M\in Z_L$ implies $M_{\left( +\right) }\in \Omega _L.$
\end{proposition}

\begin{proposition}
\label{p12} For any $L\in {\Bbb M}_n$ there exists a unique $P\in $ $%
\Psi {\bf D}_q$ of the form $P=D+\sum_{i=0}^\infty p_iD^{-i}$ such that $%
P^n=L.$
\end{proposition}

We set $P=L^{1/n}.$

\begin{proposition}
\label{p13} Any element $M\in Z_L$ is uniquely represented as 
\begin{equation*}
M=\sum_{i=-\infty }^{m(M\text{)}}c_iL^{\frac in},\quad c_i\in {\Bbb C}.
\end{equation*}
\end{proposition}

Propositions \ref{p11}, \ref{p12} imply that Lax equations 
\begin{equation}
\frac{dL}{dt}=\left[ A,L\right] ,\quad L\in {\Bbb M}_n,\quad
A=M_{\left( +\right) },\quad M=\sum_{i=-\infty }^{m(M\text{)}}c_iL^{\frac
in},\quad c_i\in {\Bbb C},  \label{f121}
\end{equation}
are self-consistent; without loss of generality we may assume that $c_i=0$
if $l\vert i.$

\begin{lemma}
\label{l14} Equation (\ref{f121}) implies that 
\begin{equation*}
\frac d{dt}L^{\frac rn}=\left[ A,L^{\frac rn}\right] \;\text{for any }r\in 
{\Bbb N}.
\end{equation*}
\end{lemma}

Lemma \ref{l14} immediately implies

\begin{proposition}
\label{p15} Functionals 
\begin{equation*}
H_m=\frac nm{\rm Tr}\;L^{\frac mn},\quad m\in {\Bbb N},
\end{equation*}
are conservation laws for (\ref{f121}).
\end{proposition}

\begin{proposition}
\label{p16} Let 
\begin{eqnarray*}
\frac{dL}{dt} &=&\left[ M_{\left( +\right) },L\right] ,\quad
M=\sum_{i=-\infty }^{m(M)}c_iL^{\frac in},\quad c_i\in {\Bbb C}, \\
\frac{dL}{d\tau } &=&\left[ \tilde M_{\left( +\right) },L\right] ,\quad
\tilde M=\sum_{i=-\infty }^{m(\tilde M)}\tilde c_iL^{\frac in},\quad \tilde
c_i\in {\Bbb C},
\end{eqnarray*}
be two Lax equations associated with any two elements in $Z_L.$ Then 
\begin{equation*}
\frac{d^2L}{dtd\tau }=\frac{d^2L}{d\tau dt};
\end{equation*}
in other words, (formal) flows generated by $M,$ $\tilde M$ commute with
each other.
\end{proposition}

\subsection{q-difference Lax equations as Hamiltonian systems}

In this section we shall describe a family $\left\{ ,\right\} _n,\;n\in 
{\Bbb N},$ of Poisson structures on $\Psi {\bf D}_q;$ the bracket $\left\{
,\right\} _n$ may be restricted to ${\Bbb M}_n\subset \Psi {\bf D}_q$
and Lax equations (\ref{lax}) are Hamiltonian with respect to this bracket.
We shall see later that $\left\{ ,\right\} _n$ coincides with the q-deformed
Gelfand-Dickey bracket \cite{FrRes, SemSev} associated with 
the Lie algebra ${\frak sl}(n).$

An accurate definition of the Poisson structure should begin with the
description of a class of admissible functionals and of their derivatives.
In the present context the algebra of observables ${\cal A}$ is generated by
'elementary' functionals which assign to a pseudodifference operator $A$ the
formal integrals of its coefficients, 
\begin{equation*}
\zeta _i^j\left( A\right) ={\rm Tr\;}\left( z^{-j}AD^{-i}\right) .
\end{equation*}

By definition, a functional $\varphi \in {\cal A}$ is smooth if for each $%
L\in {\Bbb M}_n\subset \Psi {\bf D}_q$ there exists an element $X\in
\Psi {\bf D}_q$ (called its linear gradient) such that 
\begin{equation*}
\left\langle d\varphi (L),X\right\rangle =\left( \frac d{dt}\right)
_{t=0}\varphi \left( L+tX\right) .
\end{equation*}
In applications, various functionals may be defined only on an affine
subspace of $\Psi {\bf D}_q;$ in that case the choice of the gradient (when
it exists) is not unique (however, a canonical choice is frequently
possible). It is easy to see that 'elementary' functionals are smooth; in a
similar way, traces of fractional powers of a pseudodifference operator are
smooth functionals defined on affine subspaces ${\Bbb M}_n;$ the
gradient of such a functional may be so chosen that 
\begin{equation*}
\left[ d\varphi (L),L\right] =0.
\end{equation*}
Along with the linear gradient of a functional we shall frequently use its
left and right gradients $\nabla ,\nabla ^{\prime }$ which are formally
defined by 
\begin{equation*}
\begin{array}{l}
\left\langle \nabla \varphi (L),X\right\rangle =\left( \frac d{dt}\right)
_{t=0}\varphi \left( \left( 1+tX\right) L\right) ,\\
\left\langle \nabla
^{\prime }\varphi (L),X\right\rangle =\left( \frac d{dt}\right)
_{t=0}\varphi \left( L\left( 1+tX\right) \right) ;
\end{array}
\end{equation*}
obviously, $\nabla \varphi (L)=Ld\varphi (L),\quad \nabla ^{\prime }\varphi
(L)=d\varphi (L)L.$ A functional $\varphi \in {\cal A}$ is called invariant
if $\nabla \varphi =\nabla ^{\prime }\varphi .$

Let us put ${\frak d}=\Psi {\bf D}_q\oplus \Psi {\bf D}_q$ (direct sum of
two copies); we introduce an invariant inner product in ${\frak d}$ by 
\begin{equation}
\left\langle \left\langle \left( 
\begin{array}{l}
X_1 \\ 
X_2
\end{array}
\right) \,,\left( 
\begin{array}{l}
Y_1 \\ 
Y_2
\end{array}
\right) \right\rangle \right\rangle =\left\langle X_1,Y_1\right\rangle
-\left\langle X_2,Y_2\right\rangle .  \label{f130}
\end{equation}
For a functional $\varphi $ let us write \ $D\varphi =\left( \nabla \varphi
,\nabla ^{\prime }\varphi \right) \in {\frak d.}$ We shall consider a class
of Poisson brackets on $\Psi {\bf D}_q$ which depend bilinearly on left and
right gradients of their arguments. In a very general way, such a bracket
may be written as 
\begin{equation*}
\{\varphi ,\psi \}=\left\langle \left\langle RD\varphi ,D\psi \right\rangle
\right\rangle ,
\end{equation*}
where $R\in End{\frak d}$ , \ $R=$ $\left( 
\begin{array}{cc}
A & B \\ 
C & D
\end{array}
\right)$.\footnote{Poisson brackets of this type were discussed 
by L.Freidel and J.-M.Maillet \cite{FM} and by L.Li and S.Parmentier \cite{LP}.}

We shall postpone the discussion of the Jacobi identity for this class of
brackets until part~3. Note only that it holds for all brackets constructed
below.

A natural additional condition on this class of brackets is the involutivity
of invariant functionals. It is easy to see that this condition (which
allows to use formal traces to generate commuting Hamiltonian flows) is
equivalent to the following simple constraint: 
\begin{equation*}
A+B=C+D.
\end{equation*}
A similar class of Poisson brackets is also defined in the
pseudodifferential case. In this latter case, there is a simple standard
choice of the operators $A,B,C,D$: $A=D,$ $B=C=0;$ moreover, the operators $%
A=D$ should be skew symmetric and satisfy the modified classical Yang-Baxter
equation. The standard choice is $%
A=\frac 12\left( P_{\left( +\right) }-P_{-}\right) $ (it corresponds to the
second Gelfand-Dickey bracket, which is a special case of the general
Sklyanin bracket). In the q-pseudodifference case this simple choice is no
longer possible; indeed, the standard classical r-matrix 
\begin{equation*}
r_s=\frac 12\left( P_{\left( +\right) }-P_{-}\right)
\end{equation*}
is no longer skew, because of the different properties of the invariant
inner product. Since the symmetric part of $r_s$ is the projection operator
onto the subspace of operators of order zero, it is natural to look for
modified brackets of the form 
\begin{equation}
\left\{ \varphi ,\psi \right\} =\left\langle \left\langle \left( 
\begin{array}{cc}
r+aP_0 & bP_0 \\ 
cP_0 & r+dP_0
\end{array}
\right) D\varphi ,D\psi \right\rangle \right\rangle ,  \label{f131}
\end{equation}
where $r=\frac 12\left( P_{+}-P_{-}\right) $ and $a,b,c,d$ are linear
operators acting in $J_0$ which satisfy 
\begin{equation*}
a=-a^{*},\quad d=-d^{*},\quad c^{*}=b.
\end{equation*}
In other words, the bracket (\ref{f131}) differs from the naive
Gelfand-Dickey bracket by a 'perturbation term' which is acting only on the
constant terms of the gradients (cf. \cite{S}). We shall see below that for any choice of $%
a,b,c,d$ this bracket satisfies the Jacobi identity. The additional
conditions which allow to fix the choice of the bracket completely are given
by the following uniqueness theorem.

\begin{theorem}
\label{t17}There exists a unique Poisson bracket of the form (\ref{f131}) on 
$\Psi {\bf D}_q$ such that

\begin{itemize}
\item[1)]  the affine subspace ${\Bbb M}_n$ is a Poisson submanifold;

\item[2)]  the Hamiltonians $H_m=\frac nm{\rm Tr}\;L^{\frac mn},\quad m\in {\Bbb %
N},$ are in involution and give rise to Lax equations 
\begin{equation}
\frac{dL}{dt}=\left[ L_{\left( +\right) }^{\frac mn},L\right] ,\quad L\in 
{\Bbb M}_n.  \label{f133}
\end{equation}
\end{itemize}

This bracket is given by 
\begin{equation}
\left\{ \varphi ,\psi \right\} =\left\langle \left\langle \left( 
\begin{array}{cc}
r+\frac 12\frac{1+\hat h^n}{1-\hat h^n}P_0^{\prime } & -\frac{\hat h^n}{%
1-\hat h^n}P_0^{\prime }+\frac 12P_{00} \\ 
\frac 1{1-\hat h^n}P_0^{\prime }+\frac 12P_{00} & r-\frac 12\frac{1+\hat h^n%
}{1-\hat h^n}P_0^{\prime }
\end{array}
\right) D\varphi ,D\psi \right\rangle \right\rangle .  \label{f134}
\end{equation}
\end{theorem}

\begin{remark}
Although the Poisson bracket satisfying the conditions of the theorem is
unique, there remains some freedom in the choice of the corresponding
r-matrix. The reason is that the gradients $D\varphi ,D\psi $ are not
arbitrary, namely, they belong to a family of isotropic linear subspaces in $%
{\frak d}$ ; hence $R$ is defined only up to an operator whose bilinear form
identically vanishes on all such subspaces. As an example note that the
bracket (\ref{f134}) may be also written in the form 
\begin{equation}
\left\{ \varphi ,\psi \right\} =\left\langle \left\langle \left( 
\begin{array}{cc}
P_{+}+\frac 1{1-\hat h^n}P_0^{\prime } & -\frac{\hat h^n}{1-\hat h^n}%
P_0^{\prime } \\ 
\frac 1{1-\hat h^n}P_0^{\prime } & P_{+}-\frac{\hat h^n}{1-\hat h^n}%
P_0^{\prime }
\end{array}
\right) D\varphi ,D\psi \right\rangle \right\rangle .  \label{f134p1}
\end{equation}
\end{remark}

\begin{remark}
\label{r12} Just as in the pseudodifferential case we may linearize the
quadratic bracket (\ref{f134}) at the unit element of $\Psi {\bf D}_q;$ the
resulting bracket $\left\{ \cdot ,\cdot \right\} _1$ is linear; it is given
by 
\begin{equation}
\left\{ \varphi ,\psi \right\} _1\left( X\right) =-\left\langle \left[
r_sd\varphi ,d\psi \right] +\left[ d\varphi ,r_sd\psi \right]
,X\right\rangle ,  \label{f134p2}
\end{equation}
i.e., it is the Lie-Poisson bracket associated with the r-matrix $r_s.$ The
brackets (\ref{f134}) and (\ref{f134p2}) are compatible, i.e., their linear
combinations are also Poisson brackets. Thus we have a 1-parameter family of
quadratic Poisson brackets: 
\begin{equation}
\left\{ \varphi ,\psi \right\} _\alpha =\left\{ \varphi ,\psi \right\}
+\alpha \left\{ \varphi ,\psi \right\} _1.  \label{f134p3}
\end{equation}

As usual, dynamical systems generated by the Hamiltonians $H_m$ are {\em %
bihamiltonian}; namely, the vector field generated by $H_m$ with respect to
the quadratic bracket (\ref{f134}) coincides with the vector field generated
by $H_{m+n}$ with respect to the linear bracket (\ref{f134p2}). Functionals $%
H_m,$ $m\leq n$, are Casimir functions for the bracket (\ref{f134p2}).
\end{remark}

{\em Proof of the theorem}. The gradients of $H_m$ are given by 
\begin{equation*}
\nabla H_m=\nabla ^{\prime }H_m=L^{\frac mn},
\end{equation*}
hence the Hamiltonian equation generated by $H_m$ with respect to the
bracket (\ref{f131}) is given by 
\begin{equation*}
\frac{dL}{dt}=\left( \left[ r+\left( a+b\right) P_0\right] \,L^{m/n}\right)
\cdot L-L\cdot \left( \left[ r+\left( c+d\right) P_0\right] \,L^{m/n}\right)
;
\end{equation*}
since $\left[ L,L^{\frac mn}\right] =0,$ we get 
\begin{equation*}
\begin{array}{rcl}
\dsize\frac{dL}{dt} & = &\left( \left[ P_{\left( +\right) }+\left( a+b-\frac 12\right)
P_0\right] \,L^{m/n}\right) \cdot L\\
& &\quad -L\cdot \left( \left[ P_{\left( +\right)
}+\left( c+d-\frac 12\right) P_0\right] \,L^{m/n}\right) .
\end{array}
\end{equation*}
This equation reduces to the Lax form (\ref{f133}) if and only if the
coefficients $a,b,c,d$ are such that for any $L\in {\Bbb M}_n$%
\begin{equation}
\left( \left[ a+b-\frac 12\right] \,\left( L^{m/n}\right) _0\right) \cdot
L=L\cdot \left( \left[ c+d-\frac 12\right] \,\left( L^{m/n}\right) _0\right)
.  \label{identity}
\end{equation}

\begin{lemma}
\label{l19} Condition (\ref{identity}) implies that $a+b-1/2=\left(
c+d-1/2\right) =F,$ where $F$ is a linear operator in $J_0$ with $%
ImF\subseteq $ ${\Bbb C}\cdot 1\subset J_0.$
\end{lemma}

Lemma \ref{l19} together with the antisymmetry condition imply that 
\begin{equation}
b=\frac 12-a+F,\quad c=\frac 12+a+F^{*},\quad d=-a+F-F^{*}.  \label{f142}
\end{equation}
It is easy to see that the bilinear form of the operators $F,F^{*}$ vanishes
on the gradients $D\varphi =\left( \nabla \varphi ,\nabla ^{\prime }\varphi
\right) $ of arbitrary functionals and hence does not contribute to the
Poisson bracket; indeed, 
$$
\begin{array}{rcl}
\left[ 
\begin{array}{l}
\text{the contribution from} \\ 
F,F^{*} \text{\ to}\quad \left\{ \varphi ,\psi \right\}
\end{array}
\right] &=&\left\langle \left\langle \left( 
\begin{array}{cc}
0 & FP_0 \\ 
F^{*}P_0 & \left( F-F^{*}\right) P_0
\end{array}
\right) D\varphi ,D\psi \right\rangle \right\rangle \\
&=&FP_0\left( \nabla ^{\prime }\varphi \right) \cdot \left( {\rm Tr}\nabla
\psi -{\rm Tr}\nabla ^{\prime }\psi \right)+\\ 
& &\quad +FP_0\left( \nabla ^{\prime
}\psi \right) \cdot \left( {\rm Tr}\nabla \varphi -{\rm Tr}\nabla ^{\prime
}\varphi \right) \\
&=&0,\text{ due to invariance of Tr.}
\end{array}
$$
Thus we get 
\begin{equation}
\left\{ \varphi ,\psi \right\} =\left\langle \left\langle \left( 
\begin{array}{cc}
P_{+}+\left( \frac 12+a\right) P_0 & \left( \frac 12-a\right) P_0 \\ 
\left( \frac 12+a\right) P_0 & P_{+}+\left( \frac 12-a\right) P_0
\end{array}
\right) D\varphi ,D\psi \right\rangle \right\rangle .  \label{f143}
\end{equation}
The condition that the affine subspace $\ {\Bbb M}_n$ is a Poisson
subvariety allows to fix the remaining free operator $a.$ This condition
means that the functionals 
\begin{equation}
\varphi _f\left( L\right) =\int \frac{dz}zu_n\left( z\right) f\left(
z\right) \equiv {\rm Tr}\left( LD^{-n}f\right) ,\quad \forall f\in {\Bbb C}%
((z^{-1}))  \label{f144}
\end{equation}
are Casimir functions on ${\Bbb M}_n$, i.e., 
\begin{equation}
\{\varphi _f,\psi \}\mid _{{\Bbb M}_n}=0  \label{f145}
\end{equation}
for any $\psi \in {\cal A}.$ From (\ref{f143}) we get 
\begin{equation}
\{\varphi _f,\psi \}=\left\langle \left[ \frac 12+a\right] \left( \nabla
\varphi _f\right) _0+\left[ \frac 12-a\right] \left( \nabla ^{\prime
}\varphi _f\right) _0\,,\left( \nabla \psi \right) _0-\left( \nabla ^{\prime
}\psi \right) _0\right\rangle .  \label{f146}
\end{equation}
Note that for any $L\in {\Bbb M}_n$ 
\begin{equation*}
\left( \nabla \varphi _f\right) _0=f,\quad \left( \nabla ^{\prime }\varphi
_f\right) _0=\hat h^{-n}\left( f\right) ,
\end{equation*}
in other words, the constant terms of the gradients do not depend on $L$. 
Hence the condition (\ref{f145}) is reduced to
the following one:
\begin{itemize}
\item For any $f\in {\Bbb C}((z^{-1}))$ and any $L\in {\Bbb M}_n$  
\begin{equation}
\left\langle \left[ \frac 12+a\right] f
+\left[ \frac 12-a\right] \; {}^{h^{-n}}f,
\left( \nabla \psi \right) _0-\left( \nabla ^{\prime
}\psi \right) _0\right\rangle =0 .  \label{f146p}
\end{equation}
\end{itemize}
The latter condition is reduced to
\begin{equation}
\left[ \frac 12\left( 1+\hat h^{-n}\right) +a\left( 1-\hat h^{-n}\right)
\right] f\in {\Bbb C}\quad \text{for all }f\in 
{\Bbb C}((z^{-1})). \label{f146p1} 
\end{equation} Indeed, (\ref{f146p1}) results immediately from the following  
\begin{lemma}
For any $g\in {\Bbb C}((z^{-1}))$ such that $\int \frac{dz}zg(z)=0$ there
exists a functional $\psi _g\in {\cal A}$ such that for some $L\in {\bf 
{\Bbb M}_n}$ 
\begin{equation*}
\left( \nabla \psi \left( L\right) \right) _0-\left( \nabla ^{\prime }\psi
\left( L\right) \right) _0=g. \qquad \qquad \bigtriangleup
\end{equation*} 
\end{lemma}

\noindent The proof of this assertion is similar to that of lemma \ref{l311} below.

\noindent Conversely,  (\ref{f146p1}) implies (\ref{f146p}) due to invariance of the 
trace.

From (\ref{f146p1}) we obtain that $a=a_0+\beta -\gamma ^{*},$ where 
\begin{equation*}
a_0=\frac 12\frac{1+h^n}{1-h^n}\left( 1-P_{00}\right)
\end{equation*}
and $\beta ,\gamma $ are one-dimensional linear operators in ${\Bbb C}%
((z^{-1}))\ $with $Im\beta ,\gamma \subset {\Bbb C}\cdot 1.$ The
antisymmetry of $a$ implies $\beta =\gamma .$ Thus we have 
\begin{equation*}
\begin{array}{ll}
a=a_0+\gamma -\gamma ^{*}, & b=\frac 12-a_0+F-\gamma +\gamma ^{*}, \\ 
c=\frac 12+a_0+F^{*}-\gamma ^{*}+\gamma ,\quad & d=-a_0+F-\gamma
-F^{*}+\gamma ^{*}.
\end{array}
\end{equation*}
To conclude the proof let us observe that $F,\gamma $ do not contribute to
the Poisson bracket, which implies (\ref{f134}). More precisely, we have the
following assertion:

\begin{lemma}
\label{l110p} Let $f,g,h,k$ be linear operators in $J_0$ with images in the
subspace of constants ${\Bbb C}\cdot 1\subset J_0.$ The r-matrices $R$ and $%
R^{\prime }=R+\Delta $ where 
\begin{equation*}
\Delta =\left( 
\begin{array}{cc}
h-k^{*} & f+k^{*} \\ 
h+g^{*} & -g^{*}+f
\end{array}
\right) ,
\end{equation*}
define the same Poisson bracket.
\end{lemma}

\vspace{1cm}

Now we shall prove that the bracket (\ref{f134}) coincides with the
q-deformed Gelfand-Dickey bracket derived in \cite{FRS, SemSev} via the
q-deformed Drinfeld-Sokolov reduction procedure. Let us first of all briefly
recall this reduction procedure.

Let us denote by $L{\frak gl}(n)$ the loop algebra associated with ${\frak gl%
}(n),$ i.e., the algebra of $n\times n$ matrices with coefficients in ${\Bbb %
C}((z^{-1})).$

It is well known that a scalar q-difference equation of order $n$ 
\begin{equation*}
L\psi _0=0,\quad L=D^n+u_{n-1}(z)D^{n-1}+\cdot \cdot \cdot +u_0(z),
\end{equation*}
is equivalent to a first order matrix equation 
\begin{equation*}
D\Psi ={\cal L}\Psi ,\quad \Psi =\left( 
\begin{array}{c}
\psi _0 \\ 
\vdots \\ 
\psi _{n-1}
\end{array}
\right) ,
\end{equation*}
where the potential ${\cal L\!\!}\in L{\frak gl}(n)$ has a special form. The
standard choice for ${\cal L}$ is given by a companion matrix, 
\begin{equation}
{\cal L}=\left( 
\begin{array}{cccc}
0 & 1 & \cdots & 0 \\ 
\vdots & \vdots & \ddots & \vdots \\ 
0 & 0 & \cdots & 1 \\ 
-u_0 & -u_1 & \cdots & -u_{n-1}
\end{array}
\right) .  \label{f152}
\end{equation}
This choice is not unique; a linear change of variables 
\begin{equation*}
\Psi \mapsto \Psi ^{\prime }=S\Psi ,
\end{equation*}
where $S$ a lower triangular matrix with coefficients in ${\Bbb C}((z^{-1}))$%
, induces a gauge transformation 
\begin{equation}
{\cal L\longmapsto L}^{\prime }=^hS{\cal L\,}S^{-1}.  \label{f154}
\end{equation}
Let us denote by ${\Bbb Y}_n\subset L{\frak gl}(n)$ the subvariety of all
matrices of the form 
\begin{equation}
{\cal L}^{\prime }=\left( 
\begin{array}{cccc}
\ast & 1 & \cdots & 0 \\ 
\vdots & \vdots & \ddots & \vdots \\ 
\ast & * & \cdots & 1 \\ 
\ast & * & \cdots & *
\end{array}
\right) .  \label{f155}
\end{equation}
It is easy to see that the gauge action (\ref{f154}) of the group $L{\bf N}%
_{-}\left( n\right) $ of lower triangular matrices with the coefficients in $%
{\Bbb C}((z^{-1}))$ preserves ${\Bbb Y}_n$.

\begin{theorem}
\label{orbits}\cite{FRS,SemSev}

\begin{enumerate}
\item  The gauge action of $L{\bf N}_{-}\left( n\right) $ on ${\Bbb Y}_n$ is
free.

\item  The set of companion matrices of the form (\ref{f152}) is a
cross-section of this action.
\end{enumerate}
\end{theorem}

This theorem implies that the quotient ${\Bbb Y}_{n}/L{\bf N}_{-}\left(
n\right) $ can be identified with ${\Bbb M}_{n}$.

In \cite{FRS,SemSev} a natural description of the quotient ${\Bbb Y}_n/L{\bf %
N}_{-}\left( n\right) $ in the framework of Poisson reduction has been
proposed. Let us recall some basic notions.

Let ${\cal M}$ be a Poisson manifold. The action of a Lie group $G$ on $%
{\cal M}$ is called {\it admissible} if the ring of G-invariant functions $%
I_G\left( {\cal M}\right) $ is a Poisson subalgebra in $Fun\left( {\cal M}%
\right) $ . Assume that the quotient ${\cal M}/G$ is a smooth manifold, then 
$I_G\left( {\cal M}\right) \approx Fun\left( {\cal M}/G\right) $ and
therefore the quotient ${\cal M}/G$ has a Poisson structure. The natural
projection $\pi :{\cal M\rightarrow M}/G$ is Poisson with respect to this
bracket.

\begin{proposition}
\label{constraints}Let $V\subset {\cal M}$ be a submanifold preserved by the
action of $G$. The quotient $V/G$ is a Poisson submanifold in ${\cal M}/G$
if and only if the ideal $I_0\subset I_G\left( {\cal M}\right) $ of all
G-invariant functions vanishing on $V$ is a Poisson ideal in $I_G\left( 
{\cal M}\right) $.
\end{proposition}

In our setting ${\cal M}=L{\frak gl}(n)$, $G=L{\bf N}_{-}\left( n\right) ,$ $%
V={\Bbb Y}_n$ . In order to define the q-deformed Drinfeld-Sokolov reduction
we need to find a Poisson structure on $L{\frak gl}(n)$ which satisfies the
following conditions:

\begin{enumerate}
\item  the gauge action of $L{\bf N}_{-}\left( n\right) $ on $L{\frak gl}(n)$
is admissible;

\item  the constraints defining the submanifold ${\Bbb Y}_n\subset L{\frak gl%
}(n)$ generate a Poisson ideal in $I_{L{\bf N}_{-}\left( n\right) }\left( L%
{\frak gl}(n)\right) .$
\end{enumerate}

\noindent The latter condition means that the Poisson brackets of the
constraints with any function vanish on the constraints surface ${\Bbb Y}%
_n\subset L{\frak gl}(n)$, i.e., the constraints are of the {\em first
class, }according to Dirac.

As shown in \cite{FRS, SemSev}, these two conditions allow to fix the
Poisson structure on $L{\frak gl}(n)$ and the underlying classical r- matrix
in an essentially unique way. To give an explicit formula for this bracket
let us fix the following notation.

Let $L{\frak n}_{+}\left( n\right) ,L{\frak n}_{-}\left( n\right) ,L{\frak h}%
\left( n\right) $ be the subalgebras of strictly upper triangular, strictly
lower triangular and diagonal matrices in $L{\frak gl}(n)$ respectively; let 
${\cal P}_{+},{\cal P}_{-,}{\cal P}_0$ be the corresponding projectors. Let $%
R_s$ be the automorphism of $L{\frak h}\left( n\right) $ given by 
\begin{equation*}
R_s{\bf diag}\left( \alpha _0,\ldots ,\alpha _{n-1}\right) ={\bf diag}\left(
\alpha _{n-1},\alpha _0,\ldots ,\alpha _{n-2}\right) ;
\end{equation*}
(this is the automorphism of $L{\frak h}$ induced by the {\em Coxeter element%
} of the Weyl group). Put $\theta =R_s\hat h.$ The r-matrix

\begin{equation}
\hat r=\frac 12\left( {\cal P}_{+}-{\cal P}_{-}+\frac{1+\theta }{1-\theta }%
{\cal P}_0^{\prime }\right) ,  \label{f156}
\end{equation}
where 
\begin{equation}
{\cal P}_0^{\prime }={\cal P}_0-\frac 1n\int \frac{dz}z{\rm Tr},
\label{f157}
\end{equation}
satisfies modified classical Yang-Baxter equation and is skew symmetric with
respect to the invariant inner product 
\begin{equation}
\left\langle A\left( z\right) ,B\left( z\right) \right\rangle =\int \frac{dz}%
z{\rm Tr}A\left( z\right) B\left( z\right) ,\quad A\left( z\right) ,B\left(
z\right) \in L{\frak gl}(n).  \label{f158}
\end{equation}

The Poisson bracket on $L{\frak gl}(n)$ which makes possible the q-deformed
Drinfeld-Sokolov reduction is given by 
\begin{equation}
\left\{ \hat{\varphi},\hat{\psi}\right\} _{S}\text{=}\left\langle
\left\langle \left( 
\begin{array}{cc}
\hat{r} & -\hat{h}\hat{r}_{+} \\ 
\hat{r}_{-}\hat{h}^{-1} & -\hat{r}
\end{array}
\right) \left( 
\begin{array}{l}
\nabla \hat{\varphi} \\ 
\nabla ^{\prime }\hat{\varphi}
\end{array}
\right) ,\left( 
\begin{array}{l}
\nabla \hat{\psi} \\ 
\nabla ^{\prime }\hat{\psi}
\end{array}
\right) \right\rangle \right\rangle  \label{f159}
\end{equation}
where $\hat{r}_{\pm }=\hat{r}\pm \frac{1}{2}id.$

On the reduced space ${\Bbb Y}_n/L{\bf N}_{-}\left( n\right) $ which we
identify with ${\Bbb M}_n$ we obtain a Poisson bracket called the q-deformed
(second) Gelfand-Dickey structure.

\begin{theorem}
\label{t112}The bracket (\ref{f134}) coincides with the q-deformed
Gelfand-Dickey structure $\left\{ \cdot ,\cdot \right\} _q.$
\end{theorem}

{\em Proof}. For a functional $f$ on $L{\frak gl}(n)$ put 
\begin{equation}
Z_f=^{h^{-1}}\nabla f-\nabla ^{\prime }f.  \label{f161}
\end{equation}
Let us denote by $L{\frak b}_{-}\left( n\right) $ the subalgebra of lower
triangular matrices in $L{\frak gl}(n)$ with arbitrary diagonal elements.

\begin{lemma}
\label{l113}A functional $f$ is $L{\bf N}_{-}\left( n\right) $ - invariant
if and only if $Z_f\in L{\frak b}_{-}\left( n\right) .$
\end{lemma}

\begin{lemma}
\label{l114}The value of $\left\{ f,g\right\} _S\left( {\cal L}\right) $ at
any ${\cal L}\in {\Bbb Y}_n$ does not depend on the $L{\frak n}_{-}\left(
n\right) $-components of $df,dg$ provided that $Z_f,Z_g\in L{\frak b}%
_{-}\left( n\right) .$
\end{lemma}

For any $L\in {\Bbb M}_n$ we denote by ${\cal L}$ the corresponding
companion matrix of the form (\ref{f152}). Let $\varphi ,\psi $ be any
functionals on ${\Bbb M}_n;$ by construction, the quotient Poisson structure
on ${\Bbb Y}_n/L{\bf N}_{-}\left( n\right) \simeq $ ${\Bbb M}_n$ is given by 
\begin{equation}
\left\{ \varphi ,\psi \right\} _q\left( L\right) =\left\{ \hat \varphi ,\hat
\psi \right\} _S\left( {\cal L}\right),  \label{f160}
\end{equation}
where $\hat \varphi ,\hat \psi $ are any $L{\bf N}_{-}\left( n\right) $ -
invariant functionals on $L{\frak gl}(n)$ such that their restrictions on $%
{\Bbb Y}_n$ coincide with the pullbacks of $\varphi $ and $\psi ,$
respectively: 
\begin{equation}
\hat \varphi \mid _{{\Bbb Y}_n}=\pi ^{*}\varphi ,\quad \hat \psi \mid _{%
{\Bbb Y}_n}=\pi ^{*}\psi .  \label{f160a}
\end{equation}
To calculate the r.h.s. of (\ref{f160}) we need to know only the gradients $%
d\hat \varphi ,$ $d\hat \psi $ and not $\hat \varphi ,$ $\hat \psi $
themselves. The upper triangular components of the gradients are fixed by (%
\ref{f160a}), and their strictly lower triangular components may be chosen
arbitrarily, provided that $Z_{\hat \varphi }\,,Z_{\hat \psi }\in $ $L{\frak %
b}_{-}\left( n\right) ,$ in agreement with lemma \ref{l114}.

Note that $\varphi $, $\psi $ are defined only on ${\Bbb M}_n,$ hence their
gradients are defined modulo the annihilator $\widehat{{\Bbb M}_{\,}}_{n-1}$%
of the tangent space to ${\Bbb M}_n$ . To fix them we shall suppose that
they have the form 
\begin{equation}
d\varphi =\sum_{i=0}^{n-1}f_iD^{-i},\quad f_i\in {\Bbb C}((z^{-1})).
\label{f162}
\end{equation}

\begin{lemma}
\label{l115}The upper triangular component of $d\hat \varphi $ is given by 
\begin{equation}
d\hat \varphi _{pm}\left( {\cal L}\right) =-{\rm Res}\left( D^pd\varphi
\left[ LD^{-(m+1)}\right] _{(+)}\right) ,\quad m\geq p.  \label{f163}
\end{equation}
\end{lemma}

Let us define the strictly lower triangular components of $d\hat \varphi $
by the same formula. We need to verify that $Z_{\hat \varphi }$ $\in L{\frak %
b}_{-}\left( n\right) ;$ this results directly from the following lemma:

\begin{lemma}
\label{l116} We have 
\begin{equation}
\label{f164}
\begin{array}{rcl}
\nabla \hat \varphi _{pm}\left( {\cal L}\right) &=&\delta _{n-1,p}\left(
\nabla \varphi \left[ LD^{-(m+1)}\right] _{(+)}\right) _0-\\
& & \quad\quad -\bar \delta
_{n-1,p}\left( D^{p+1}d\varphi \left[ LD^{-(m+1)}\right] _{(+)}\right) _0;
\end{array}
\end{equation}
\begin{eqnarray}
\left( Z_{\hat \varphi }\right) _{n-1,0}\left( {\cal L}\right)
&=&^{h^{-1}}\left( \nabla \varphi \left[ LD^{-1)}\right] _{(+)}\right)
_0-\left( D^{n-1}d\varphi \right) _0u_0;  \label{f165} \\
\left( Z_{\hat \varphi }\right) _{n-1,m}\left( {\cal L}\right)
&=&^{h^{-1}}\left( \nabla \varphi \left[ LD^{-(m+1)}\right] _{(+)}\right)
_0,\quad m\neq 0;  \label{f166} \\
\left( Z_{\hat \varphi }\right) _{p-1,0}\left( {\cal L}\right) &=&-\left(
D^p\nabla ^{\prime }\varphi \right) _0,  \label{f167}
\end{eqnarray}
where $\bar \delta =1-\delta .$ All other elements of $Z_{\hat \varphi }$
are zero.
\end{lemma}

Taking into account that $Z_{\hat \varphi }$ $,Z_{\hat \psi }\in L{\frak b}%
_{-}\left( n\right) ,$ we get the following expression for the bracket (\ref
{f159}) 
\begin{equation}
\left\{ \hat \varphi ,\hat \psi \right\} _S=\frac 12\left( \left\langle 
\frac{1+\theta }{1-\theta }Z_{\hat \varphi }^0,Z_{\hat \psi }^0\right\rangle
+\left\langle ^hZ_{\hat \varphi },\nabla \hat \psi \right\rangle
-\left\langle \nabla \hat \varphi ,^hZ_{\hat \psi }\right\rangle \right) ,
\label{f168}
\end{equation}
where $Z_{\hat \varphi }^0\equiv {\cal P}_0^{\prime }Z_{\hat \varphi },\quad
Z_{\hat \psi }^0\equiv {\cal P}_0^{\prime }Z_{\hat \psi }.$

Let us calculate the contribution of the first term in (\ref{f168}) to $%
\left\{ \hat \varphi ,\hat \psi \right\} _S$ .

\begin{lemma}
\label{l117}The eigenfunctions of the operator $\theta $ are 
\begin{equation}
E_{m,\alpha }=z^m{\bf e}_\alpha ,\quad m\in {\Bbb Z}\;,\quad \alpha
=0,\ldots ,n-1,  \label{f169}
\end{equation}
where 
\begin{equation}
{\bf e}_\alpha ={\bf diag}\left( 1,\omega ^{-\alpha },\ldots ,\omega
^{-\left( n-1\right) \alpha }\right) ,\quad \omega =e^{\frac{2\pi i}n}.
\label{f170}
\end{equation}
The corresponding eigenvalues $\xi _{m,\alpha }$ are equal to 
\begin{equation}
\xi _{m,\alpha }=q^m\omega ^\alpha .  \label{f171}
\end{equation}
The eigenfunctions satisfy the condition 
\begin{equation}
\left\langle E_{m,\alpha },E_{l,\beta }\right\rangle =n\delta _{m,-l}\cdot
\left\{ 
\begin{array}{l}
1,\quad \alpha =-\beta \ {\rm \limfunc{mod}\;}n,{\rm \ } \\ 
0,\quad \text{in the other cases},
\end{array}
\right.  \label{f172}
\end{equation}
and form a basis in $L{\frak h}\left( n\right) .$
\end{lemma}

We shall denote by $\sum_{m,\alpha }^{\prime }$ the sum over all pairs $%
\left( m,\alpha \right) \neq \left( 0,0\right) ,\quad m\in {\Bbb Z}\;,\quad
\alpha =0,\ldots ,n-1.$ Note that in the expansion of $Z_{\hat \varphi
}^0,Z_{\hat \psi }^0$ with respect to the eigenbasis $E_{m,\alpha }$ the $%
E_{0,0}$-component is absent, hence 
\begin{equation*}
\left\langle \frac{1+\theta }{1-\theta }Z_{\hat \varphi }^0,Z_{\hat \psi
}^0\right\rangle =\left\langle \mathop{{\sum}'}\limits_{m,\alpha }\frac 1n\frac{%
1+q^m\omega ^\alpha }{1-q^m\omega ^\alpha }E_{m,\alpha }\left\langle Z_{\hat
\varphi }\,,E_{-m,n-\alpha }\right\rangle \;,\,Z_{\hat \psi }\right\rangle=
\end{equation*}
\begin{equation}
=\int \frac{dz}z\int \frac{dw}w\mathop{{\sum}'}\limits_{m,\alpha }\frac 1n\frac{%
1+q^m\omega ^\alpha }{1-q^m\omega ^\alpha }\left( \frac zw\right) ^m{\rm Tr}%
\left( {\cal P}_0Z_{\hat \varphi }\left( w\right) \cdot {\bf e}_{n-\alpha
}\right) {\rm Tr}\left( {\cal P}_0Z_{\hat \psi }\left( z\right) \cdot {\bf e}%
_\alpha \right).  \label{f173}
\end{equation}
Applying lemma \ref{l116} we obtain: 
\begin{equation}
\label{f174}
\begin{array}{l}
{\rm Tr}\left( {\cal P}_0Z_{\hat \varphi }\left( w\right) \cdot {\bf e}%
_{n-\alpha }\right) {\rm Tr}\left( {\cal P}_0Z_{\hat \psi }\left( z\right)
\cdot {\bf e}_\alpha \right) =   \\
\quad =^{h^{-1}}P_0\nabla \varphi \left( w\right) \cdot \;^{h^{-1}}P_0\nabla \psi
\left( z\right) +P_0\nabla ^{\prime }\varphi \left( w\right) \cdot P_0\nabla
^{\prime }\psi \left( z\right) -   \\
\quad \quad -\omega ^{-\alpha }\cdot \;^{h^{-1}}P_0\nabla \varphi \left( w\right) \cdot
P_0\nabla ^{\prime }\psi \left( z\right) -\omega ^\alpha P_0\nabla ^{\prime
}\varphi \left( w\right) \cdot \;^{h^{-1}}P_0\nabla \psi \left( z\right) .
\end{array}
\end{equation}
We denote by $A_1^1,\ldots ,A_1^4$ the contributions of the corresponding
terms of (\ref{f174}) to $\left\langle \frac{1+\theta }{1-\theta }Z_{\hat
\varphi }^0,Z_{\hat \psi }^0\right\rangle .$

\begin{lemma}
\label{l118} 
\begin{eqnarray}
\frac 1n\sum_{\alpha =0}^{n-1}\frac{1+q^m\omega ^\alpha }{1-q^m\omega
^\alpha } &=&\frac{1+q^{mn}}{1-q^{mn}},\quad m\neq 0;  \label{f175} \\
\frac 1n\sum_{\alpha =0}^{n-1}\frac{1+\omega ^\alpha }{1-\omega ^\alpha }
&=&0  \label{f176} \\
\frac 1n\sum_{\alpha =0}^{n-1}\frac{1+q^m\omega ^\alpha }{1-q^m\omega
^\alpha }\omega ^\alpha &=&2\frac{q^{m\left( n-1\right) }}{1-q^{mn}},\quad
m\neq 0;  \label{f177} \\
\frac 1n\sum_{\alpha =0}^{n-1}\frac{1+\omega ^\alpha }{1-\omega ^\alpha }%
\omega ^\alpha &=&-\frac{n-2}n  \label{f178} \\
\frac 1n\sum_{\alpha =0}^{n-1}\frac{1+q^m\omega ^\alpha }{1-q^m\omega
^\alpha }\omega ^{-\alpha } &=&2\frac{q^m}{1-q^{mn}},\quad m\neq 0;
\label{f179} \\
\frac 1n\sum_{\alpha =0}^{n-1}\frac{1+\omega ^\alpha }{1-\omega ^\alpha }%
\omega ^{-\alpha } &=&\frac{n-2}n  \label{f180}
\end{eqnarray}
\end{lemma}

{\em Proof}. Let us prove (\ref{f175}); formulae (\ref{f177}), (\ref{f179})
may be verified in the same way. We have 
\begin{equation*}
S_1=\frac 1n\sum_{\alpha =0}^{n-1}\frac{1+q^m\omega ^\alpha }{1-q^m\omega
^\alpha }=-1+\frac 2n\sum_{\alpha =0}^{n-1}\frac 1{1-q^m\omega ^\alpha
}=-1+\frac 2n\sum_{i=0}^\infty \sum_{\alpha =0}^{n-1}q^{mi}\omega ^{\alpha
i},
\end{equation*}
but 
\begin{equation*}
\frac 1n\sum_{\alpha =0}^{n-1}\omega ^{\alpha i}=\left\{ 
\begin{array}{l}
0,\quad i\neq jn,\quad j\in {\Bbb N}, \\ 
1,\quad i=jn,
\end{array}
\right.
\end{equation*}
and hence 
\begin{equation*}
S_1=-1+2\sum_{j=0}^\infty q^{mnj}=\frac{1+q^{mn}}{1-q^{mn}},
\end{equation*}
as desired.

To prove (\ref{f176}) note that 
\begin{equation*}
\frac{1+\omega ^{n-\alpha }}{1-\omega ^{n-\alpha }}=\frac{1+\omega ^{-\alpha
}}{1-\omega ^{-\alpha }}=-\frac{1+\omega ^\alpha }{1-\omega ^\alpha },
\end{equation*}
therefore for odd $n$ all terms in the sum (\ref{f176}) cancel each other
completely. If $n$ is even only the term $\dfrac{1+\omega ^{n/2}}{1-\omega
^{n/2}}$ survives, but evidently it is zero, since $\omega ^{n/2}=-1$.
Formulae (\ref{f178}), (\ref{f180}) immediately follow from (\ref{f176}).

\begin{lemma}
\label{l118p} 
\begin{eqnarray}
A_1^1 &=&\left\langle \frac{1+\hat h^n}{1-\hat h^n}P_0^{\prime }\nabla
\varphi ,\,\nabla \psi \right\rangle ,  \label{f181} \\
A_1^2 &=&\left\langle \frac{1+\hat h^n}{1-\hat h^n}P_0^{\prime }\nabla
^{\prime }\varphi ,\,\nabla ^{\prime }\psi \right\rangle ,  \label{f182} \\
A_1^3 &=&\left\langle \frac{-2}{1-\hat h^n}P_0^{\prime }\nabla \varphi
,\,\nabla ^{\prime }\psi \right\rangle -\frac{n-2}n{\rm Tr}\nabla \varphi
\cdot {\rm Tr}\nabla ^{\prime }\psi ,  \label{f183} \\
A_1^4 &=&\left\langle \frac{-2\hat h^n}{1-\hat h^n}P_0^{\prime }\nabla
^{\prime }\varphi ,\,\nabla \psi \right\rangle +\frac{n-2}n{\rm Tr}\nabla
^{\prime }\varphi \cdot {\rm Tr}\nabla \psi .  \label{f184}
\end{eqnarray}
\end{lemma}

{\em Proof. }We verify only (\ref{f183}), other formulae can be proved in
the same way. We have: 

\begin{eqnarray*}
A_1^3 &=&-\int \frac{dz}z\int \frac{dw}w\mathop{{\sum}'}\limits_{m,\alpha }\frac 1n%
\frac{1+q^m\omega ^\alpha }{1-q^m\omega ^\alpha }\omega ^{-\alpha }\left(
\frac zw\right) ^m\cdot \; ^{h^{-1}}P_0\nabla \varphi \left( w\right) \cdot
P_0\nabla ^{\prime }\psi \left( z\right)  \\
&=&-\int \frac{dz}z\left[ \int \frac{dw}w\sum_{{m\in {\Bbb Z}}\atop{m\neq 0}%
}\frac 1n\frac{2q^m}{1-q^{mn}}\omega ^{-\alpha }\left( \frac zw\right)
^m\cdot \; ^{h^{-1}}P_0\nabla \varphi \left( w\right) \right] \cdot P_0\nabla
^{\prime }\psi \left( z\right)  \\
&&-\int \frac{dz}z\left[ \int \frac{dw}w\left( \frac 1n\sum_{\alpha =0}^{n-1}%
\frac{1+\omega ^\alpha }{1-\omega ^\alpha }\omega ^{-\alpha }\right) \cdot \;
^{h^{-1}}P_0\nabla \varphi \left( w\right) \right] \cdot P_0\nabla ^{\prime
}\psi \left( z\right) \\
&=&-\int \frac{dz}z\left( \frac{2\hat h}{1-\hat h^n}\,^{h^{-1}}P_0^{\prime
}\nabla \varphi \right) \left( z\right) \cdot P_0\nabla ^{\prime }\psi
\left( z\right) - \\
&&-\frac{n-2}n\int \frac{dw}w\,^{h^{-1}}P_0^{\prime }\nabla \varphi \left(
w\right)\cdot\int\frac{dz}z P_0\nabla ^{\prime }\psi
\left( z\right)\qquad \left( \text{by lemma \ref{l118}}\right) \\
&=&\left\langle \frac{-2}{1-\hat h^n}P_0^{\prime }\nabla \varphi ,\,\nabla
^{\prime }\psi \right\rangle -\frac{n-2}n{\rm Tr}\nabla \varphi \cdot {\rm Tr%
}\nabla ^{\prime }\psi ,
\end{eqnarray*}
as desired.

Using lemma \ref{l118p} and taking into account that ${\rm Tr}\nabla \varphi =%
{\rm Tr}\nabla ^{\prime }\varphi ,$ we obtain 
\begin{eqnarray}
\left\langle \frac{1}{2}\frac{1+\theta }{1-\theta }Z_{\hat{\varphi}}^{0},Z_{%
\hat{\psi}}^{0}\right\rangle &=&\left\langle \left\langle \left( 
\begin{array}{cc}
\frac{1}{2}\frac{1+\hat{h}^{n}}{1-\hat{h}^{n}}P_{0}^{\prime } & -\frac{\hat{h%
}^{n}}{1-\hat{h}^{n}}P_{0}^{\prime } \\ 
\frac{1}{1-\hat{h}^{n}}P_{0}^{\prime } & -\frac{1}{2}\frac{1+\hat{h}^{n}}{1+%
\hat{h}^{n}}P_{0}^{\prime }
\end{array}
\right) \left( 
\begin{array}{l}
\nabla \varphi \\ 
\nabla ^{\prime }\varphi
\end{array}
\right) ,\left( 
\begin{array}{l}
\nabla \psi \\ 
\nabla ^{\prime }\psi
\end{array}
\right) \right\rangle \right\rangle.  \notag\\
& & \label{f185}
\end{eqnarray}
It remains to calculate $\left\langle ^{h}Z_{\hat{\varphi}},\nabla \hat{\psi}%
\right\rangle -\left\langle \nabla \hat{\varphi},^{h}Z_{\hat{\psi}%
}\right\rangle $ .

\begin{lemma}
\label{l119} 
\begin{eqnarray}
\left\langle ^hZ_{\hat \varphi },\nabla \hat \psi \right\rangle
&=&\left\langle r\nabla \varphi ,\nabla \psi \right\rangle -\left\langle
r\nabla ^{\prime }\varphi ,\nabla ^{\prime }\psi \right\rangle +  \notag \\
&&\quad+\frac 12{\rm Tr}\left( P_0\nabla \varphi \cdot P_0\nabla \psi
\right) -\frac 12{\rm Tr}\left( P_0\nabla ^{\prime }\varphi \cdot P_0\nabla
^{\prime }\psi \right).  \label{f186}
\end{eqnarray}
\end{lemma}

{\em Proof}. Taking into account that $\left[ LD^{-n}\right] _{\left(
+\right) }=1$ and using lemma \ref{l116}, we obtain 
\begin{equation}
\left\langle ^{h}Z_{\hat{\varphi}},\nabla \hat{\psi}\right\rangle
=A_{2}^{1}+A_{2}^{2}+A_{2}^{3}+A_{2}^{4},  \label{f187}
\end{equation}
where 
\begin{eqnarray*}
A_{2}^{1} &=&{\rm Tr}\left( \sum_{p=0}^{n-2}\left( D^{p+1}\nabla ^{\prime
}\varphi D^{-1}\right) _{0}\left( Dd\psi \left[ LD^{-\left( p+1\right)
}\right] _{\left( +\right) }\right) _{0}\right) , \\
A_{2}^{2} &=&-{\rm Tr}\left( \sum_{p=0}^{n-2}\left( \nabla \varphi \left[
LD^{-\left( m+1\right) }\right] _{\left( +\right) }\right) _{0}\left(
D^{m+1}d\psi \right) _{0}\right) , \\
A_{2}^{3} &=&{\rm Tr}\left( ^{h}\left[ \left( D^{n-1}d\varphi \right)
_{0}u_{0}\right] \left( Dd\psi \right) _{0}\right) , \\
A_{2}^{4} &=&{\rm Tr}\left( \left( \nabla \varphi \right) _{0}\left( \nabla
\psi \right) _{0}\right) .
\end{eqnarray*}

In transformations below we use the fact that ${\rm Tr}A_{\left( +\right) }B=%
{\rm Tr}AB_{\left( -\right) }$ and the following proposition:

\begin{proposition}
\label{p119p}Let $B\in \Psi {\bf D}_q$ be of the form 
\begin{equation*}
B=\sum_{i=m}^lb_iD^{-i},\quad b_i\in {\Bbb C}\left( \left( z^{-1}\right)
\right) ,
\end{equation*}
then 
\begin{equation*}
\sum_{i=m}^lD^{-i}\left( D^iB\right) _0=B.
\end{equation*}
\end{proposition}

Since $d\varphi ,$ $d\psi $ have the form (\ref{f162}), we obtain: 
\begin{eqnarray}
A_2^1 &=&{\rm Tr}\left( Dd\psi \left[ L\sum_{p=0}^{n-2}D^{-\left( p+1\right)
}\left( D^{p+1}\nabla ^{\prime }\varphi D^{-1}\right) _0\right] _{\left(
+\right) }\right)  \notag \\
&=&{\rm Tr}\left( Dd\psi \left[ L\left( \nabla ^{\prime }\varphi
D^{-1}\right) _{-}-LD^{-n}\left( D^n\nabla ^{\prime }\varphi D^{-1}\right)
_0\right] _{\left( +\right) }\right)  \notag \\
&=&{\rm Tr}\left( \left( Dd\psi \right) _{\left( -\right) }L\left( \nabla
^{\prime }\varphi D^{-1}\right) _{-}\right) -{\rm Tr}\left( Dd\psi \left(
D^n\nabla ^{\prime }\varphi D^{-1}\right) _0\right)  \notag \\
&=&{\rm Tr}\left( \left( Dd\psi -D\left( d\psi \right) _0\right) L\left(
\nabla ^{\prime }\varphi D^{-1}\right) _{-}\right) -{\rm Tr}\left( \left(
Dd\psi \right) _0\cdot \; ^h\left( D^{n-1}\nabla ^{\prime }\varphi \right)
_0\right)  \notag \\
&=&{\rm Tr}\left( \left( \nabla ^{\prime }\psi -\left( d\psi \right)
_0L\right) \nabla ^{\prime }\varphi _{\left( -\right) }\right) -{\rm Tr}%
\left( \left( Dd\psi \right) _0\cdot \; ^h\left[ \left( D^{n-1}d\varphi \right)
_{\left( -\right) }L\right] \right)  \notag \\
&=&{\rm Tr}\nabla ^{\prime }\varphi \nabla ^{\prime }\psi _{\left( +\right)
}-{\rm Tr}\left( d\psi \right) _0L\nabla ^{\prime }\varphi -{\rm Tr}\left(
\left( Dd\psi \right) _0\cdot \; ^h\left[ \left( D^{n-1}d\varphi \right)
_{\left( 0\right) }u_0\right] \right)  \notag \\
&=&{\rm Tr}\nabla ^{\prime }\varphi \nabla ^{\prime }\psi _{\left( +\right)
}-{\rm Tr}\nabla \varphi L\left( d\psi \right) _0-A_2^3.  \label{f188}
\end{eqnarray}
$A_2^2$ may be developed as follows: 
\begin{eqnarray}
A_2^2 &=&-{\rm Tr}\left( \sum_{p=0}^{n-2}\left( \nabla \varphi \left[
LD^{-\left( m+1\right) }\right] _{\left( +\right) }\right) _0\left(
D^{m+1}d\psi \right) _0\right)  \notag \\
&=&-{\rm Tr}\left( \nabla \varphi \left[ \sum_{p=0}^{n-2}LD^{-\left(
m+1\right) }\left( D^{m+1}d\psi \right) _0\right] _{\left( +\right) }\right)
\notag \\
&=&-{\rm Tr}\left( \nabla \varphi \left[ L\left( d\psi -\left( d\psi \right)
_0\right) \right] _{\left( +\right) }\right)  \notag \\
&=&-{\rm Tr}\nabla \varphi \left( \nabla \psi \right) _{\left( +\right) }+%
{\rm Tr}\nabla \varphi L\left( d\psi \right) _0 . \label{f189}
\end{eqnarray}
Substituting (\ref{f188}), (\ref{f189}) in (\ref{f187}) gives 
\begin{equation*}
\left\langle ^hZ_{\hat \varphi },\nabla \hat \psi \right\rangle ={\rm Tr}%
\nabla ^{\prime }\varphi \nabla ^{\prime }\psi _{\left( +\right) }-{\rm Tr}%
\nabla \varphi \left( \nabla \psi \right) _{\left( +\right) }+{\rm Tr}\left(
\nabla \varphi \right) _0\left( \nabla \psi \right) _0
\end{equation*}
which immediately implies (\ref{f186}). Lemma \ref{l119} is proved.

For $\left\langle \nabla \hat \varphi ,^hZ_{\hat \psi }\right\rangle $ we
have a relation similar to (\ref{f186}), hence 
\begin{equation}
\left\langle ^hZ_{\hat \varphi },\nabla \hat \psi \right\rangle
-\left\langle \nabla \hat \varphi ,^hZ_{\hat \psi }\right\rangle =2\left(
\left\langle r\nabla \varphi ,\nabla \psi \right\rangle -\left\langle
r\nabla ^{\prime }\varphi ,\nabla ^{\prime }\psi \right\rangle \right) .
\label{f190}
\end{equation}
Substituting (\ref{f190}) and (\ref{f185}) in (\ref{f168}) we obtain 
\begin{equation}
\left\{ \varphi ,\psi \right\} _q=\left\langle \left\langle \left( 
\begin{array}{cc}
r+\frac 12\frac{1+\hat h^n}{1-\hat h^n}P_0^{\prime } & -\frac{\hat h^n}{%
1-\hat h^n}P_0^{\prime } \\ 
\frac 1{1-\hat h^n}P_0^{\prime } & r-\frac 12\frac{1+\hat h^n}{1+\hat h^n}%
P_0^{\prime }
\end{array}
\right) \left( 
\begin{array}{l}
\nabla \varphi \\ 
\nabla ^{\prime }\varphi
\end{array}
\right) ,\left( 
\begin{array}{l}
\nabla \psi \\ 
\nabla ^{\prime }\psi
\end{array}
\right) \right\rangle \right\rangle.  \label{f191}
\end{equation}
This formula differs from (\ref{f134}) by the absence of $\pm \frac 12P_{00}$
in the non-diagonal elements of r-matrix, but the corresponding terms give
no contribution to the bracket because of the invariance of the inner
product: 
\begin{equation*}
\left\langle \frac 12P_{00}\nabla ^{\prime }\varphi ,\nabla \psi
\right\rangle -\left\langle \frac 12P_{00}\nabla \varphi ,\nabla ^{\prime
}\psi \right\rangle =\frac 12{\rm Tr}\nabla ^{\prime }\varphi {\rm Tr}\nabla
\psi -{\rm Tr}\nabla \varphi {\rm Tr}\nabla ^{\prime }\psi =0.
\end{equation*}
Theorem \ref{t112} is proved.

\begin{theorem}
\label{t120}The Poisson bracket (\ref{f134}) in terms of generating
functions 
\begin{equation}
u_i\left( z\right) =\sum_{m=-\infty }^{N\left( u_i\right) }u_{im}z^m,\quad
i=0,\ldots ,n-1,  \label{f192p}
\end{equation}
has the form 
\begin{eqnarray}
\left\{ u_i\left( z\right) ,u_j\left( w\right) \right\} &=&\sum_{{m\in 
{\Bbb Z}}\atop{m\neq 0}}\frac{\left( 1-q^{m\left( n-i\right) }\right) \left(
1-q^{mj}\right) }{1-q^{mn}}\left( \frac wz\right) ^mu_i\left( z\right)
u_j\left( w\right) +  \notag \\
&& +\sum_{r=1}^{\min \left( n-i,j\right) }\delta \left( \frac{wq^r}%
z\right) u_{i+r}\left( w\right) u_{j-r}\left( z\right) -  \notag \\
&& -\sum_{r=1}^{\min \left( n-i,j\right) }\delta \left( \frac
w{zq^{i-j+r}}\right) u_{i+r}\left( z\right) u_{j-r}\left( w\right) ,
\label{f192}
\end{eqnarray}
where $\delta \left( z\right) =\sum\limits_{m\in {\Bbb Z}}z^m.$ So this
bracket coincides with the one constructed by Frenkel and Reshetikhin in 
\cite{FrRes}.
\end{theorem}

Proof of this theorem is straightforward computation.

\section{The group of q-pseudodifference symbols of all complex degrees and
the associated q-KdV hierarchies.}

\subsection{The double extension of the algebra $\Psi {\bf D}_q$ .}

Let us define the operator $\ln D\in {\rm End}\left( {\Bbb C}\left( \left(
z^{-1}\right) \right) \right) $ by 
\begin{equation}
\ln D=\ln q\cdot z\frac d{dz},  \label{f31}
\end{equation}
where the branch of $\ln q$ is fixed by 
\begin{equation}
-\pi <\arg q<\pi ,\quad \ln 1=0.  \label{f31p}
\end{equation}
As above, we suppose that $|q|<1$ . Note that the subspaces ${\Bbb C}z^m$
are the eigenspaces for $\ln D$ with eigenvalues $\lambda _m=m\ln q$ , hence
the exponential $\exp \ln D$ is well-defined and $\exp \ln D=D,$ which
justifies our definition. Evidently, 
\begin{equation}
\left[ \ln D,f\right] =\ln q\cdot z\frac{df}{dz},\quad \forall f\in {\Bbb C}%
\left( \left( z^{-1}\right) \right) ,\qquad \left[ \ln D,D\right] =0,
\label{f32}
\end{equation}
which implies that $\left[ \ln D,\cdot \right] $ is an (outer) derivation of
the associative algebra $\Psi {\bf D}_q.$

\begin{proposition}
\label{p31}The 2-form 
\begin{equation}
\omega \left( X,Y\right) =\left\langle \left[ \ln D,X\right] ,Y\right\rangle
,\quad \forall X,Y\in \Psi {\bf D}_q,  \label{f33}
\end{equation}
is a 2-cocycle on $\Psi {\bf D}_q.$
\end{proposition}

{\em Proof}. Note that for any $X\in \Psi {\bf D}_q$ 
\begin{equation*}
\limfunc{Tr}\left[ \ln D,X\right] =0,
\end{equation*}
which, together with (\ref{f32}), implies the skew-symmetry of $\omega .$
Next we have 
\begin{equation*}
\left[ \ln D,\left[ X,Y\right] \right] =\left[ X,\left[ \ln D,Y\right]
\right] -\left[ Y,\left[ \ln D,X\right] \right] ,
\end{equation*}
hence

\begin{eqnarray*}
\omega \left( \left[ X,Y\right] ,Z\right) &=&\left\langle \left[ X,\left[
\ln D,Y\right] \right] ,Z\right\rangle -\left\langle \left[ Y,\left[ \ln
D,X\right] \right] ,Z\right\rangle \\
&=&\left\langle \left[ \ln D,Y\right] ,\left[ Z,X\right] \right\rangle
+\left\langle \left[ \ln D,X\right] ,\left[ Y,Z\right] \right\rangle \\
&=&-\omega \left( \left[ Z,X\right] ,Y\right) -\omega \left( \left[
Y,Z\right] ,X\right) ,
\end{eqnarray*}
as desired.

The logarithmic cocycle $\omega $ defines a non-trivial central extension $%
\widehat{\Psi {\bf D}_q}=\Psi {\bf D}_q\dotplus {\Bbb C\cdot }{\bf c}$ of
the Lie algebra $\Psi {\bf D}_q.$ $\widehat{\Psi {\bf D}_q}$ does not admit
a non-degenerate invariant inner product. To improve the situation let us
consider the ''double extension'' $\widetilde{\Psi {\bf D}_q}$ of the
algebra $\Psi {\bf D}_q:$%
\begin{equation}
\widetilde{\Psi {\bf D}_q}=\Psi {\bf D}_q\dotplus {\Bbb C\cdot }\ln
D\dotplus {\Bbb C\cdot } {\bf c};  \label{f38}
\end{equation}
the bilinear form 
\begin{equation}
\left\langle X+\alpha \ln D+\beta {\bf c},Y+\gamma \ln D+\delta {\bf c}%
\right\rangle =\left\langle X,Y\right\rangle _{\Psi {\bf D}_q}+\alpha \delta
+\beta \gamma  \label{f39}
\end{equation}
is invariant, non-degenerate and sets into duality the subspaces $J_{+}$ and 
$J_{-};$ moreover, 
\begin{equation}
\quad J_0^{*}\approx J_0,\quad \left( {\Bbb C\cdot }\ln D\right) ^{*}\approx 
{\Bbb C\cdot }{\bf c}  \label{f310}
\end{equation}
(here $J_{\pm },J_0$ are defined by (\ref{f15})--(\ref{f17}) ).

\subsection{The group of q-pseudodifference symbols of all complex degrees.}

For any $\alpha \in {\Bbb C}$ we define the complex power $\hat h^\alpha $
of the operator $\hat h$ by 
\begin{equation}
\left( \hat h^\alpha f\right) \left( z\right) \equiv \,^{h^\alpha }f\left(
z\right) =f\left( q^\alpha z\right) ,  \label{f313}
\end{equation}
where $q^\alpha =\exp \left( \alpha \ln q\right) $ and the branch of $\ln q$
is fixed by (\ref{f31p}).

A normalized q-pseudodifference symbol of degree $\alpha $ is a formal
series of the form 
\begin{equation}
L=D^\alpha +\sum_{i=1}^\infty a_iD^{\alpha -i},\quad a_i\in {\Bbb C}\left(
\left( z^{-1}\right) \right) .  \label{f314}
\end{equation}
The multiplication law of symbols is uniquely defined by the commutation
relation: 
\begin{equation}
D^\alpha \circ a=\,^{h^\alpha }a\circ D^\alpha .  \label{f315}
\end{equation}
Let $\widehat{G}_\alpha $ be the set of normalized q-difference symbols of
degree $\alpha $ and $\widehat{G}_{-}$ be the set of symbols of all complex
degrees, 
\begin{equation*}
\widehat{G}_{-}=\bigcup\limits_{\alpha \in {\Bbb C}}\widehat{G}_\alpha .
\end{equation*}
$\widehat{G}_{-}$ is a group with respect to the multiplication law (\ref
{f315}). This group admits the following description. For $\alpha \in {\Bbb C%
}$ let $\sigma _\alpha $ be the automorphism of $\widehat{G}_{-}$ given by 
\begin{equation}
\sigma _\alpha \left( X\right) =D^\alpha XD^{-\alpha }.  \label{f314p}
\end{equation}
Obviously, $\sigma _\alpha $ preserves the degree of symbols and hence
induces an automorphism of the subgroup $\widehat{G}_0.$

\begin{lemma}
\label{l31p}The group $\widehat{G}_{-}$ is the semi-direct product of the
additive group ${\Bbb C}$ and the group $\widehat{G}_0$.
\end{lemma}

\noindent For $L\in \widehat{G}_{-}$ we shall write $L=\bar LD^\alpha ,$ where $%
\alpha $ is the degree of $L$ and $\bar L\in $ $\widehat{G}_0$.

\noindent In an obvious sense, $\widehat{G}_{-}$ may be regarded as an
infinite-dimensional Lie group.

\begin{lemma}
\label{l32}The tangent Lie algebra of the group $\widehat{G}_{-}$ is the
algebra $\hat J_{-}=J_{-}\dotplus {\Bbb C\cdot }\ln D$ considered as a Lie
subalgebra in $\widetilde{\Psi {\bf D}_q}.$
\end{lemma}

{\em Proof.} We must check that 
\begin{equation}
\left[ \ln D,X\right] =\frac d{d\alpha }\mid _{\alpha =0}D^\alpha
XD^{-\alpha }  \label{f318p}
\end{equation}
for all $X\in \Psi {\bf D}_q.$ Since $\left[ D^\alpha ,D\right] =0,$ it is
sufficient to prove (\ref{f318p}) for $X=f,\quad \forall f\in {\Bbb C}\left(
\left( z^{-1}\right) \right) .$ In this case we have 
\begin{equation*}
\frac d{d\alpha }\mid _{\alpha =0}D^\alpha f\left( z\right) D^{-\alpha
}=\frac d{d\alpha }\mid _{\alpha =0}f\left( q^\alpha z\right) =\ln q\cdot z%
\frac{df}{dz}=\left[ \ln D,f\right] ,
\end{equation*}
as desired. $\blacksquare$

Let us fix the following models of the tangent and cotangent spaces of $%
\widehat{G}_{-}$ which will be used over the rest of this section: 
\begin{equation}
\begin{array}{lll}
T_L\widehat{G}_{-} & = & \left\{ X=\overline{X}D^\alpha +\widetilde{X}L\ln
D,\quad \overline{X}\in J_{-},\quad \widetilde{X}\in {\Bbb C}\right\} , \\ 
T_L^{*}\widehat{G}_{-} & = & \left\{ f=D^{-\alpha }\bar f+\tilde f{\bf c}%
L^{-1},\quad \bar f\in J_{+},\quad \tilde f\in {\Bbb C}\right\} .
\end{array}
\label{f319}
\end{equation}
The pairs $\left( \overline{X},\widetilde{X}\right) $ and $\left( \bar
f,\tilde f\right) $ will be called {\em normal coordinates }of $X$ and $f,$
respectively. The canonical pairing between $T_L\widehat{G}_{-}$ and $T_L^{*}%
\widehat{G}_{-}$ is given by 
\begin{equation}
\left\langle X,f\right\rangle =\left\langle \overline{X},\bar f\right\rangle
_{\Psi {\bf D}_q}+\widetilde{X}\cdot \tilde f.  \label{f320p}
\end{equation}
Let us denote by $^lS_L$ and $^rS_L$ the operators of the left (resp.,
right) multiplication by $L$ in $\widehat{G}_{-}.$

\begin{lemma}
\label{l33}In normal coordinates the tangent map $\left( ^rS_L\right) _{*}$
at the unit element of $\widehat{G}_{-}$%
\begin{equation*}
\left( ^rS_L\right) _{*\left( e\right) }:T_e\widehat{G}_{-}\rightarrow T_L%
\widehat{G}_{-},\quad X\mapsto Y,
\end{equation*}
is given by: 
\begin{equation}
\begin{array}{lll}
\overline{Y} & = & \overline{X}\bar L+\widetilde{X}\left[ \ln D,\bar
L\right] , \\ 
\widetilde{Y} & = & \widetilde{X}.
\end{array}
\label{f322}
\end{equation}
The tangent map $\left( ^lS_L\right) _{*\left( e\right) }$ is given by 
\begin{equation}
\begin{array}{lll}
\overline{Y} & = & \bar L\sigma _\alpha \left( \overline{X}\right) , \\ 
\widetilde{Y} & = & \widetilde{X}.
\end{array}
\label{f323}
\end{equation}
Formulae (\ref{f323}) may be formally written as $Y=LX$.
\end{lemma}

\begin{proposition}
\label{p34}
 \hspace*{1cm} \linebreak

1) The exponential map $\exp :\hat J_{-}\rightarrow \widehat{G}_{-}$ is well
defined on the whole $\hat J_{-}.$

2) The restriction of the exponential map to the affine subspace $\hat
J_\alpha =J_{-}+\alpha \ln D$ is a bijection between $\hat J_\alpha $ and $%
\widehat{G}_\alpha $ if 
\begin{equation}
\frac{\alpha \ln q}{2\pi i}\notin {\Bbb Q}.  \label{f324}
\end{equation}
\end{proposition}

{\em Proof}. By definition, $L\left( t\right) =\exp t\left( X+\ln D\right)
,\quad X\in J_{-},$ is a solution of the following differential equation 
\begin{equation}
\frac{dL}{dt}=\left( ^lS_L\right) _{*\left( e\right) }\left( X+\ln D\right)
\label{f326}
\end{equation}
with the initial condition $L\left( 0\right) =e.$ Evidently, $L\left(
t\right) $ has the form $L\left( t\right) =A\left( t\right) D^t,$ where $%
A\left( t\right) \in \widehat{G}_0.$ By lemma \ref{l33} we write (\ref{f326}%
) as 
\begin{equation}
\frac{dA}{dt}=A\sigma _t\left( X\right) ,\quad t\in {\Bbb C}.  \label{f327}
\end{equation}
This equation has a unique solution with initial condition $A\left( 0\right)
=e.$ Indeed, let $A_i^j\left( t\right) ,\;X_i^j$ be the coefficients of the
expansion of $A\left( t\right) $ (respectively, $X$ ) in terms of $z$ and $D$
: 
\begin{equation*}
\begin{array}{rllrll}
A\left( t\right) & = & \sum\limits_{i=0}^\infty A_i\left( t\right)
D^{-i},\quad & A_i\left( t\right) & = & \sum\limits_{j=-\infty
}^{m_i}A_i^j\left( t\right) z^j,\;j\geq 1,\;A_0=1; \\ 
X & = & \sum\limits_{i=1}^\infty X_iD^{-i},\quad & X_i & = & 
\sum\limits_{j=-\infty }^{n_i}X_i^jz^j.
\end{array}
\end{equation*}
(We set $A_i^j\left( t\right) \equiv 0,\;j>m_i,\quad X_i^j\equiv 0,\;j>n_i,$
and extend summation up to infinity.) From (\ref{f327}) we obtain 
\begin{equation}
\begin{array}{lll}
\dfrac{dA_1^m}{dt} & = & q^{tm}X_1^m,\quad A_1^m\left( 0\right) =0;
\end{array}
\label{f330}
\end{equation}
\begin{equation}
\begin{array}{lll}
\dfrac{dA_i^m}{dt} & = & \sum\limits_{j=o}^{i-1}\sum\limits_{n\in {\Bbb Z}%
}A_j^nq^{nt}X_{i-j}^{m-n}+q^{mt}X_i^m,\quad A_i^m\left( 0\right) =0.
\end{array}
\label{f332}
\end{equation}

We shall prove inductively that this system admits a unique solution which
is holomorphic in ${\Bbb C}$ . Observe first of all that the function $%
w\longmapsto q^{mw}$ is holomorphic in ${\Bbb C}$, hence the value of the
integral $\int_0^tq^{mw}dw$ does not depend on the path of integration and 
\begin{equation}
A_1^m\left( t\right) =\int\limits_0^tq^{mw}dw\cdot X_1^m  \label{f331}
\end{equation}
gives the unique solution of (\ref{f330}). Obviously, $A_1^m\left( t\right) $
is holomorphic in ${\Bbb C}.$

Assume now that all coefficients $A_1,\ldots ,A_{i-1}$ are holomorphic
functions. We shall deduce from it that the equation (\ref{f332}) for $A_i$
is also solvable in holomorphic functions. Indeed, the sum over $n$ in the
r.h.s. of (\ref{f332}) has only a finite number of non-zero terms and the
functions $A_j^n$ are holomorphic by the inductive hypothesis, hence the
r.h.s. of (\ref{f332}) is holomorphic and this equation has the unique
solution which is given by 
\begin{equation}
A_i^m\left( t\right) =\int\nolimits_0^t\left(
\sum\limits_{j=o}^{i-1}\sum\limits_{n\in {\Bbb Z}}A_j^n\left( w\right)
q^{nw}X_{i-j}^{m-n}+q^{mw}X_i^m\right) dw.  \label{f333}
\end{equation}
Thus, the exponential map is well-defined. To verify the second assertion of
the proposition we need to prove that for any $L\in \widehat{G}_\alpha ,$ $%
\frac{\alpha \ln q}{2\pi i}\notin {\Bbb Q},$ there is a unique
representation of the form 
\begin{equation*}
L\equiv AD^\alpha =\exp \alpha \left( X+\ln D\right) ,\quad X\in J_{-}.
\end{equation*}
Consider (\ref{f331}), (\ref{f333}) as an equation for $X.$ We have 
\begin{equation}
A_1^m=\int_0^\alpha q^{mw}dw\cdot X_1^m=\left\{ 
\begin{array}{l}
\frac{e^{\alpha m\ln q}-1}{m\ln q}X_1^m,\quad m\neq 0, \\[0.2cm] 
\alpha X_1^0,\quad m=0.
\end{array}
\right.  \label{f334}
\end{equation}
Since$\frac{\alpha \ln q}{2\pi i}\notin {\Bbb Q},,$ $e^{\alpha m\ln q}\neq
1, $ and hence (\ref{f334}) is solvable for any $A_1^m.$ Let us assume now
that $X_1,\ldots ,X_{i-1}$ are already determined. In that case equation (%
\ref{f333}) for $X_i$ has a unique solution. Indeed, we have 
\begin{equation}
A_i^m-\int\nolimits_0^\alpha \left( \sum\limits_{j=o}^{i-1}\sum\limits_{n\in 
{\Bbb Z}}A_j^n\left( w\right) q^{nw}X_{i-j}^{m-n}\right)
dw=\int_0^\alpha q^{mw}dw\cdot X_i^m.  \label{f335}
\end{equation}
Functions $A_j^n\left( w\right) $ are defined by the formulae (\ref{f333})
and may be expressed in terms of $X_1,\ldots ,X_{i-1},$ hence the l.h.s. of (%
\ref{f335}) is a known number. As above, we see that (\ref{f335}) is
uniquely solvable provided that $\frac{\alpha \ln q}{2\pi i}\notin {\Bbb Q}.$

\begin{definition}
We call $\alpha \in {\Bbb C}$ generic if $\frac{\alpha \ln q}{2\pi i}\notin 
{\Bbb Q}$ ; we call $L\in \widehat{G}_{-}$ a generic element if its degree $%
\alpha \equiv \deg L$ is generic.
\end{definition}

\subsection{The generalized q-deformed Gelfand-Dickey structure on $\widehat{%
G}_{-}$ and related q-KdV hierarchies.}

In this section we consider Lax equations of the form 
\begin{equation}
\frac{dL}{dt}=\left[ L_{\left( +\right) }^{\frac m\alpha },L\right] ,\quad
L\in \widehat{G}_\alpha ,\quad \alpha \text{ is generic.}  \label{f336}
\end{equation}
We shall see that in some natural class of quadratic Poisson brackets on $%
\widehat{G}_\alpha $ there exists a unique one which is consistent with the
equations (\ref{f336}) (i.e., the latter are Hamiltonian with respect to
this bracket with the Hamiltonians $H_m=\frac \alpha m\limfunc{Tr}L^{\frac
m\alpha }$). For any positive integer $\alpha $ this bracket may be
restricted to ${\Bbb M}_\alpha $ and coincides there with the q-deformed
Gelfand-Dickey structure considered in part 1.

The brackets referred to are smooth with respect to the parameter $\alpha $
outside the line $\frac{2\pi i}{\ln q}{\Bbb R},$ hence we can glue them up
to a smooth Poisson structure on $\widehat{G}_{-}^{\prime
}=\dsize\bigcup\limits_{\alpha \notin \frac{2\pi i}{\ln q}{\Bbb R}}\widehat{G}%
_\alpha .$ This bracket is called the {\em generalized q-deformed (second)
Gelfand-Dickey structure}. It is uniquely defined by the following
conditions:

\begin{itemize}
\item[1)]  $\widehat{G}_\alpha $ are Poisson submanifolds;

\item[2)]  the restriction of this bracket to $\widehat{G}_\alpha $ coincides
with the unique bracket on $\widehat{G}_\alpha $ consistent with equation (%
\ref{f336}).
\end{itemize}

\noindent Let us now turn back to equation (\ref{f336}).

\begin{theorem}
\label{t35}
\end{theorem}

\begin{enumerate}
\em
\item  The equation (\ref{f336}) is self-consistent, i.e., its r.h.s. is
well-defined and belongs to the tangent space $T_L\widehat{G}_\alpha ;$

\item  The flows corresponding to different $m$ 's commute with each other;

\item  The functionals 
\begin{equation}
H_n=\frac \alpha n\limfunc{Tr}L^{\frac n\alpha },\quad n\in {\Bbb N},
\label{f337}
\end{equation}
are conservation laws for (\ref{f336}).
\end{enumerate}

{\em Proof.} $L\in \widehat{G}_\alpha $ and $\alpha $ is generic, hence in
agreement with proposition \ref{p34} there exists an $X\in J_{-}$ such that $%
L=\exp \alpha \left( X+\ln D\right) .$ For any $\beta \in {\Bbb C}$ we
define $L^\beta $ by $L^\beta =\exp \alpha \beta \left( X+\ln D\right) ;$
clearly, $L^\beta \in \widehat{G}_{\alpha \beta }.$ In particular, $L^{\frac
m\alpha }\in \widehat{G}_m,$ hence it contains only integer powers of $D$
and may be considered as an element of $\Psi {\bf D}_q.$ So the expressions $%
L_{\left( +\right) }^{\frac m\alpha }$ and $\limfunc{Tr}L^{\frac n\alpha }$
are well-defined. Note that $\left[ L^{\frac m\alpha },L\right] =0,$
therefore $\left[ L_{\left( +\right) }^{\frac m\alpha },L\right] =-\left[
L_{-}^{\frac m\alpha },L\right] $ has the form $\sum_{i=1}^\infty
a_iD^{\alpha -i},$ i.e., it is a tangent vector from $T_L\widehat{G}_\alpha
. $ The last two assertions of the theorem may be proved in the same way as
in the case of positive integer $\alpha .$ Another proof will be given below
in the frameworks of the Hamiltonian formalism.

\begin{remark}
\label{rkhlr}A similar class of equations has been constructed in \cite{KhLR}%
. The authors start with the algebra $\Psi DO_q$ of q--difference symbols of
the form 
\begin{equation}
A=\sum_{i=-\infty }^{n\left( A\right) }u_i\left( z\right) D_q^i,\quad u_i\in 
{\Bbb C}\left[ z,z^{-1}\right] ,  \label{fr31}
\end{equation}
with the commutation relation 
\begin{equation*}
D_q\circ u=\frac{u\left( qz\right) -u\left( z\right) }{q-1}+^hu\cdot D_q,
\end{equation*}
which corresponds to the definition of $D_q$ as a q-difference analogue of
the derivative $\frac d{dz}$ : 
\begin{equation*}
\left( D_qf\right) \left( z\right) =\frac{f\left( qz\right) -f\left(
z\right) }{q-1}.
\end{equation*}
Then they define the residue of a symbol $A\in \Psi DO_q$ by $\widetilde{%
{\rm Res}}A=u_{-1}\left( z\right) $ (see (\ref{fr31}) ); the trace is given
by 
\begin{equation*}
\widetilde{{\rm Tr}}A=\int \frac{dz}z\widetilde{{\rm Res}}\left( \frac
A{\left( q-1\right) D_q+1}\right) .
\end{equation*}
Let us extend $\Psi DO_q$ by assuming that the coefficients in (\ref{fr31})
are formal power series in $z^{-1}.$ It is easy to see that after this
extension $\Psi DO_q$ becomes isomorphic to $\Psi {\bf D}_q$ (as an
associative algebra) and, moreover, the definitions of the traces in $\Psi
DO_q$ and $\Psi {\bf D}_q$ coincide up to a scalar factor.

Then the authors of \cite{KhLR} construct a double extension of $\Psi DO_q$
by adjoining to it the logarithm $\log D_q$ and the corresponding 2-cocycle,
and define the group $G_q$ of q-difference symbols of the form 
\begin{equation*}
F=D_q^\alpha +\sum_{k=1}^\infty u_k\left( z\right) D_q^{\alpha -1}.
\end{equation*}
We did not found an isomorphism between $G_q$ and $\widehat{G}_{-}$ ; the
question about the relation between the generalized q-KdV equations (\ref
{f336}) and the ones constructed in \cite{KhLR} remains open.
\end{remark}

We will now examine the Hamiltonian formalism for equations (\ref{f336}).
Let us consider the set $\widehat{G}_\alpha ,$ where $\alpha $ is generic.
Just as in part 1, we consider the Poisson brackets of the form 
\begin{equation}
\left\{ \varphi ,\psi \right\} =\left\langle \left\langle \left( 
\begin{array}{cc}
r+aP_0 & bP_0 \\ 
cP_0 & r+dP_0
\end{array}
\right) \left( 
\begin{array}{c}
\overline{\nabla \varphi } \\ 
\overline{\nabla ^{\prime }\varphi }
\end{array}
\right) ,\left( 
\begin{array}{c}
\overline{\nabla \psi } \\ 
\overline{\nabla ^{\prime }\psi }
\end{array}
\right) \right\rangle \right\rangle ,  \label{f340}
\end{equation}
where $a,b,c,d$ are linear operators in $J_0$ satisfying the skew-symmetry
conditions 
\begin{equation}
a=-a^{*},\quad d=-d^{*},\quad c^{*}=b;  \label{f341}
\end{equation}
the r-matrix $r$ is defined as above: $r=\frac 12\left( P_{+}-P_{-}\right) ;%
\overline{\nabla \varphi },\;\overline{\nabla ^{\prime }\varphi }$ denote
the normal coordinates of $\nabla \varphi ,\;\nabla ^{\prime }\varphi .$
Note that for a functional $\varphi \in Fun\left( \widehat{G}_\alpha \right) 
$ the component $\overline{d\varphi }$ of its linear gradient is defined up
to an arbitrary element of $J_{\left( -\right) },$ and $\widetilde{d\varphi }
$ is arbitrary. Lemma \ref{l33} implies that 
\begin{equation}
\overline{\nabla \varphi }=\overline{L}\,\overline{d\varphi },\quad 
\overline{\nabla ^{\prime }\varphi }=\sigma _{-\alpha }\left( \overline{%
d\varphi \,}\overline{L}\right) ,  \label{f338}
\end{equation}
where $\overline{L}=LD^{-\alpha }.$ Hence, $\widetilde{d\varphi }$ gives no
contribution to the bracket (\ref{f340}).

\begin{theorem}
\label{t37}There exists a unique Poisson bracket of the form (\ref{f340}) on 
$\widehat{G}_\alpha $ which satisfies the following conditions:

1) the expression (\ref{f340}) is well-defined, i.e., it does not depend on
the $J_{\left( -\right) }$-components of $\overline{d\varphi },$ $\overline{%
d\psi };$

2) the Hamiltonians $H_m$ (see (\ref{f337}) ) give rise to the Lax equations
(\ref{f336}).

This bracket is given by 
\begin{eqnarray}
\left\{ \varphi ,\psi \right\} &=&\left\langle \left\langle \left( 
\begin{array}{cc}
r+\frac 12\frac{1+\hat h^\alpha }{1-\hat h^\alpha }P_0^{\prime } & -\frac{%
\hat h^\alpha }{1-\hat h^\alpha }P_0^{\prime }+\frac 12P_{00} \\ 
\frac 1{1-\hat h^\alpha }P_0^{\prime }+\frac 12P_{00} & r-\frac 12\frac{%
1+\hat h^\alpha }{1-\hat h^\alpha }P_0^{\prime }
\end{array}
\right) \left( 
\begin{array}{c}
\overline{\nabla \varphi } \\ 
\overline{\nabla ^{\prime }\varphi }
\end{array}
\right) ,\left( 
\begin{array}{c}
\overline{\nabla \psi } \\ 
\overline{\nabla ^{\prime }\psi }
\end{array}
\right) \right\rangle \right\rangle .\notag \\
&&  \label{f343}
\end{eqnarray}
\end{theorem}

{\em Proof.} We shall seek for $a,b,c,d$ such that the bracket (\ref{f340})
satisfies the two conditions of the theorem.

\begin{lemma}
\label{l38} 
\begin{equation}
\overline{dH_m}=D^\alpha L^{\frac m\alpha -1}.  \label{f344}
\end{equation}
\end{lemma}
Lemma \ref{l38} and (\ref{f338}) imply that 
\begin{equation}
\overline{\nabla H_m}=\overline{\nabla ^{\prime }H_m}=L^{\frac m\alpha }.
\label{f357}
\end{equation}
Substituting this into (\ref{f340}), we obtain the following

\begin{proposition}
\label{p39} Condition (2) of the theorem is equivalent the following one:
for any $L\in \widehat{G}_\alpha ,$ any $m\in {\Bbb N},$ holds 
\begin{equation}
\left( \left[ a+b-\frac 12\right] \,\left( L^{\frac m\alpha }\right)
_0\right) \cdot L=L\cdot \left( \left[ c+d-\frac 12\right] \,\left( L^{\frac
m\alpha }\right) _0\right) .  \label{f359}
\end{equation}
\end{proposition}

\begin{lemma}
\label{l310} Condition (\ref{f359}) implies that 
\begin{equation}
a+b-1/2=\left( c+d-1/2\right) =F,  \label{f360}
\end{equation}
where $F$ is a linear operator in ${\Bbb C}((z^{-1})),$ $ImF=$ ${\Bbb C}%
\cdot 1\subset {\Bbb C}((z^{-1})).$
\end{lemma}

{\em Proof}. Since $|q|<1,$ the point $1$ is generic and by proposition \ref
{p34} the map $L\mapsto L^{\frac 1\alpha }$ is a bijection between $\widehat{%
G}_\alpha $ and $\widehat{G}_{1.}$Therefore, for any $f\in {\Bbb C}%
((z^{-1})),\;f\neq 0$, there exists an $L\in \widehat{G}_\alpha $ such that $%
\left( L^{\frac 1\alpha }\right) _0=f;$ moreover, relation (\ref{f334})
implies that in the expansion 
\begin{equation}
L=D^\alpha +u_1D^{\alpha -1}+\cdots  \label{f361}
\end{equation}
the coefficient $u_1$ is nonzero. Put $\tilde F=\left( a+b-1/2\right) f$ and 
$\tilde G=\left( c+d-1/2\right) f.$ We have $\tilde FL=L\tilde G,$ which
implies $\tilde F=\hat h^\alpha \tilde G$ and $u_1\tilde F=u_1\hat h^{\alpha
-1}\tilde G.$ But $u_1\neq 0,$ hence $\tilde F=\hat h^\alpha \tilde G=\hat
h^{\alpha -1}\tilde G,$ i.e. $\tilde F\in {\Bbb C},$ $\tilde G\in {\Bbb C},$
which, together with (\ref{f359}), implies that $\tilde F=\tilde G.
\quad\blacksquare$

\noindent Using the skew-symmetry conditions (\ref{f341}) we get

\begin{equation}
b=\frac 12-a+F,\quad c=\frac 12+a+F^{*},\quad d=-a+F-F^{*}.  \label{f365}
\end{equation}
Just as in the proof of theorem \ref{t17} we can verify that $F$ does not
contribute to the Poisson bracket. The invariance of the inner product now
implies that 
\begin{eqnarray}
\left\{ \varphi ,\psi \right\} =\left\langle \left\langle \left( 
\begin{array}{cc}
P_{+}+\left( \frac 12+a\right) P_0 & \left( \frac 12-a\right) P_0 \\ 
\left( \frac 12+a\right) P_0 & P_{+}+\left( \frac 12-a\right) P_0
\end{array}
\right) \left( 
\begin{array}{c}
\overline{\nabla \varphi } \\ 
\overline{\nabla ^{\prime }\varphi }
\end{array}
\right) ,\left( 
\begin{array}{c}
\overline{\nabla \psi } \\ 
\overline{\nabla ^{\prime }\psi }
\end{array}
\right) \right\rangle \right\rangle .\notag\\
&&  \label{f366}
\end{eqnarray}
Hence $\left\{ \varphi ,\psi \right\} $ does not depend on $J_{-}$%
-components of $\overline{d\varphi },\;\overline{d\psi }.$ The requirement
that it is also independent on $J_0$-components of the gradients fixes the
choice of $a.$

\begin{lemma}
\label{l312} The bracket (\ref{f366}) does not depend on $J_0$-components of $%
\overline{d\varphi },\;\overline{d\psi }$ if and only if 
\begin{equation*}
\left[ \left( \frac 12+a\right) +\left( \frac 12-a\right) \hat h^{-\alpha
}\right] f\in {\Bbb C}.\quad \text{for }\forall f\in {\Bbb C}((z^{-1})).
\end{equation*}
\end{lemma}

{\em Proof.} Let ${\overline{d\varphi } \, \mathstrut}^{\prime }$ 
be another representative
of $\overline{d\varphi },$ ${\overline{d\varphi }\, \mathstrut}^{\prime }=\overline{%
d\varphi }+f,\quad f\in J_0,$ and $\left\{ \varphi ,\psi \right\} ^{\prime }$
be the value of the Poisson bracket corresponding to ${\overline{d\varphi }%
\, \mathstrut}^{\prime }.$ We must prove that 
\begin{equation*}
\Delta =\left\{ \varphi ,\psi \right\} ^{\prime }-\left\{ \varphi ,\psi
\right\} =0.
\end{equation*}
We have 
\begin{equation}
\Delta =\left\langle \left\langle \left( 
\begin{array}{cc}
\frac 12+a & \frac 12-a \\[0.2cm] 
\frac 12+a & \frac 12-a
\end{array}
\right) \left( 
\begin{array}{c}
\left( \bar Lf\right) _0 \\[0.2cm] 
\sigma _{-\alpha }\left( f\,\bar L\right) _0
\end{array}
\right) ,\left( 
\begin{array}{c}
\overline{\nabla \psi } \\[0.2cm] 
\overline{\nabla ^{\prime }\psi }
\end{array}
\right) \right\rangle \right\rangle .  \label{f370}
\end{equation}
Since $f\in J_0,$ we have $\left( \bar Lf\right) _0=f$ and $\sigma _{-\alpha
}\left( f\,\bar L\right) _0=\hat h^{-\alpha }f,$ hence 
\begin{equation*}
\Delta =\left\langle \left[ \left( \frac 12+a\right) +\left( \frac
12-a\right) \hat h^{-\alpha }\right] f\,,\left( \overline{\nabla \psi }%
\right) _0-\left( \overline{\nabla ^{\prime }\psi }\right) _0\right\rangle .
\end{equation*}

\begin{lemma}
\label{l311} For any $g\in {\Bbb C}((z^{-1}))$ such that $\int \frac{dz}%
zg(z)=0$ there exists a functional $\psi _g\in Fun\left( \hat G_\alpha
\right) $ such that for some $L\in \hat G_\alpha $ 
\begin{equation*}
\left( \overline{\nabla \psi _g}\left( L\right) \right) _0-\left( \overline{%
\nabla ^{\prime }\psi _g}\left( L\right) \right) _0=g.
\end{equation*}
\end{lemma}

{\em Proof.} Note that either $\alpha -1$ or $\alpha -2$ are generic. If $%
\alpha -1$ is generic, we may suppose that 
\begin{equation*}
\psi _g=\limfunc{Tr}\bar LDa_g,\quad a_g\in {\Bbb C}((z^{-1})).
\end{equation*}
It is easy to see that in this case 
\begin{equation}
\left( \overline{\nabla \psi _g}\left( L\right) \right) _0-\left( \overline{%
\nabla ^{\prime }\psi _g}\left( L\right) \right) _0=\left( 1-\hat
h^{1-\alpha }\right) \left( u_1a_g\right) ,  \label{f367}
\end{equation}
where $u_1$ is the coefficient in the expansion of $\bar L$ in powers of $%
D^{-1}:$%
\begin{equation*}
\bar L=1+u_1D^{-1}+\cdots
\end{equation*}
The condition that $\alpha -1$ is generic implies that the equation $\left(
1-\hat h^{1-\alpha }\right) \left( u_1a_g\right) =g$ is solvable for any $%
g\in {\Bbb C}((z^{-1}))$ such that $\int \frac{dz}zg(z)=0.$

If $\alpha -2$ is generic we can find $\psi _g$ in the form $\psi _g=%
\limfunc{Tr}\bar LD^2a_g.$ $\blacksquare$

$\Delta $ must vanish for any $f\in {\Bbb C}((z^{-1}))$ and any $\psi \in
Fun\left( \hat G_\alpha \right) ,$ therefore from lemma \ref{l311} it
follows that 
\begin{equation*}
\left[ \left( \frac 12+a\right) +\left( \frac 12-a\right) \hat h^{-\alpha
}\right] f\in {\Bbb C}\quad \text{for }\forall f\in {\Bbb C}((z^{-1})).
\end{equation*}
End of the proof of this theorem is just like the one of theorem \ref{t17}.
$\blacksquare$

\begin{remark}
As in part 1 we can linearize the bracket (\ref{f343}) at $L=D^\alpha $ and
construct the family of compatible Poisson structures.
\end{remark}

Now we can prove assertions 2 and 3 of theorem \ref{t35}; they immediately
result from the following

\begin{proposition}
\label{p313}Functionals $H_n=\frac \alpha n\limfunc{Tr}L^{\frac n\alpha },$ $%
n\in {\Bbb N},$ are in involution with respect to the bracket (\ref{f343}).
\end{proposition}

\begin{proposition}
\label{p314}The submanifolds $\hat G_{\alpha ,n}\subset \hat G_\alpha $ of
the symbols of the form 
\begin{equation}
L=\left( 1+\sum_{i=1}^nu_iD^{-i}\right) D^\alpha  \label{f373}
\end{equation}
are Poisson submanifolds for bracket (\ref{f343}).
\end{proposition}

{\em Proof.} It is sufficient to check that the bracket of any functional of
the form 
\begin{equation}
\varphi _{f,l}=\int \frac{dz}zu_l\,f,\quad f\in {\Bbb C}((z^{-1})),\quad l>n,
\label{f374}
\end{equation}
with any functional $\psi $ vanishes on $\hat G_{\alpha ,n}.$ Clearly, $%
\overline{d\varphi _{f,l}}=D^lf,$ hence $\overline{\nabla \varphi _{f,l}}(L)$, 
$\overline{\nabla ^{\prime }\psi _{f,l}}(L)\in J_{+}$ for all $L \in 
\hat G_{\alpha ,n}$ and we have
\begin{equation*}
\left\{ \varphi _{f,l\,},\psi \right\} =\frac 12\left( \left\langle 
\overline{\nabla \varphi _{f,l}}\,,\overline{\nabla \psi }\right\rangle
-\left\langle \overline{\nabla ^{\prime }\varphi _{f,l}}\,,\overline{\nabla
^{\prime }\psi }\right\rangle \right) =0,
\end{equation*}
as desired.

\begin{proposition}
\label{p315}The coefficient $u_n\left( z\right) $ is a Casimir function on $%
\hat G_{n,n}.$
\end{proposition}

{\em Proof.} Define $\varphi _{f,n}$ by (\ref{f374}). We shall prove that $%
\left\{ \varphi _{f,n\,},\psi \right\} $ vanishes on $\hat G_{n,n}.$ We may
write the bracket (\ref{f343}) in the form

\begin{eqnarray}
&&\left\{ \varphi _{f,n\,},\psi \right\} =\notag\\
&&=\left\langle \left\langle \left( 
\begin{array}{cc}
-P_{-}+\left( a-\frac 12\right) P_0 & \left( \frac 12-a\right) P_0 \\ 
\left( a+\frac 12\right) P_0 & -P_{-}-\left( a+\frac 12\right) P_0
\end{array}
\right) \left( 
\begin{array}{c}
\overline{\nabla \varphi } \\ 
\overline{\nabla ^{\prime }\varphi }
\end{array}
\right) ,\left( 
\begin{array}{c}
\overline{\nabla \psi } \\ 
\overline{\nabla ^{\prime }\psi }
\end{array}
\right) \right\rangle \right\rangle . \notag\\
&&\label{f375}
\end{eqnarray}
Clearly, $\overline{\nabla \varphi _{f,n}}$ , $\overline{\nabla ^{\prime
}\varphi _{f,n}}\in J_{\left( +\right) },$ hence only the $J_0$-components
of the gradients give contribution to the bracket. We have 
\begin{equation*}
\left( \overline{\nabla ^{\prime }\varphi _{f,n}}\right) _0=\sigma
_{-n}\left( \left( \overline{d\varphi }\,\bar L\right) _0\right) =\sigma
_{-n}\left( D^nf\,u_nD^{-n}\right) =f\,u_n=\left( \overline{\nabla \varphi
_{f,n}}\right) _0.
\end{equation*}
Substituting this in (\ref{f375}) we obtain $\left\{ \varphi _{f,n\,},\psi
\right\} =0,$ as desired.

Note that for integer $n$ the submanifold $\hat G_{n,n}$ may be canonically
identified with the subspace ${\Bbb M}_n\subset \Psi {\bf D}_q$ considered
in part 1 of this paper. With this identification the restriction of the
bracket (\ref{f343}) on $\hat G_{n,n}$ coincides with the bracket (\ref{f134}%
). Indeed, the r-matrices are the same and we need only to check that the
definitions of the gradients are consistent with each other. We have 
\begin{equation*}
\begin{array}{l}
\overline{\nabla \varphi }=\bar L\overline{d\varphi }=\bar LD^n\cdot D^{-n}%
\overline{d\varphi }=Ld\varphi , \\ 
\overline{\nabla ^{\prime }\varphi }=\sigma _{-n}\left( \overline{d\varphi }%
\,\bar L\right) =D^{-n}\overline{d\varphi }\cdot \bar LD^n=d\varphi L,
\end{array}
\end{equation*}
as desired.

\section{Jacobi identity for the generalized q-deformed Gelfand-Dickey
structures.}

In this section we discuss the Jacobi identity for generalized q-deformed
Gelfand-Dickey brackets (\ref{f134}), (\ref{f343}). We start with the
following general pattern, due to Gelfand and Dorfman \cite{GelfDorfm}. Let $%
{\cal A}=\Omega _0\ \ $be an associative commutative algebra, ${\frak a}%
=Der{\cal A}$ the Lie algebra of its derivations; we regard ${\cal A}$ as a $%
{\frak a}$-module and define the Chevalley complex associated with ${\cal A}$
in the standard way, 
\begin{equation*}
\begin{array}{l}
\Omega _0\overset{d}{\rightarrow }\Omega _1\overset{d}{\rightarrow }\Omega _2%
\overset{d}{\rightarrow }...,\Omega _p={\cal A}\otimes \bigwedge^p{\frak a}%
^{*}, \\ 
d\alpha (X_1,...,X_{p+1})=\sum_i(-1)^iX_i\alpha (X_1,...,\hat
X_i,...,X_{p+1})+ \\ 
\sum_{i<j}(-1)^{i+j}\alpha (X_1,...,\hat X_i,...,\hat X_j,...X_{p+1}).
\end{array}
\end{equation*}
For $X\in {\frak a}$ let $i_X:\Omega _p\rightarrow \Omega _{p-1}$be the
inner derivative, 
\begin{equation*}
i_X\alpha \left( X_1,X_2,...,X_{p-1}\right) =\alpha \left(
X,X_1,X_2,...,X_{p-1}\right) ;
\end{equation*}
for $p=1$ the coupling $\left\langle X,\alpha \right\rangle =i_X\alpha $ is
a natural bilinear pairing between ${\frak a}$ and $\Omega _0.$ Let $%
L_X=d\cdot i_X+i_X\cdot d$ be the Lie derivative. By definition, a Poisson
operator is a linear operator $H\in Hom($ $\Omega _1,{\frak a});$ the
Schouten bracket of two Poisson operators $H,K$ is a trilinear mapping $%
\Omega _0\times \Omega _0\times \Omega _0\rightarrow {\frak a}$ defined by 
\begin{equation*}
\left[ \left[ H,K\right] \right] \left( \alpha _1,\alpha _2,\alpha _3\right)
=\left\langle HL_{K\alpha _1}\alpha _2,\alpha _3\right\rangle +\left\langle
KL_{H\alpha _1}\alpha _2,\alpha _3\right\rangle +c.p.
\end{equation*}
(as usual, $c.p.$ stands for cyclic permutation). The Poisson bracket
associated with $H$ is given by 
\begin{equation}
\left\{ \varphi ,\psi \right\} =\left\langle Hd\varphi ,d\psi \right\rangle ;
\label{pbr}
\end{equation}
this bracket is skew and satisfies the Jacobi identity if and only if $H$ is
skew-symmetric and its Schouten bracket with itself is zero.

Let ${\frak J}$ be an associative algebra with a non-degenerated invariant
inner product. We choose as ${\cal A}=\Omega _0$ the algebra of smooth
functionals on ${\frak J}.$ (Recall that by definition a functional $\varphi
\in {\cal A}$ is smooth if for each $L\in {\frak J}$ there exists an element 
$X\in {\frak J}$ (called its linear gradient) such that 
\begin{equation*}
\left\langle d\varphi (L),X\right\rangle =\left( \frac d{dt}\right)
_{t=0}\varphi \left( L+tX\right) .
\end{equation*}
The left and right gradients $\nabla ,\nabla ^{\prime }$ are given by $%
\nabla \varphi (L)=Ld\varphi (L),\quad \nabla ^{\prime }\varphi (L)=d\varphi
(L)L.$ For a functional $\varphi $ we write \ $D\varphi =\left( 
\begin{array}{l}
\nabla \varphi \\ 
\nabla ^{\prime }\varphi
\end{array}
\right) \in {\frak J}\tsize\bigoplus {\frak J}.$ ) Let us define the invariant
inner product in ${\frak J}\tsize\bigoplus {\frak J}$ by 
\begin{equation}
\left\langle \left\langle \left( 
\begin{array}{l}
X_1 \\ 
X_2
\end{array}
\right) \,,\left( 
\begin{array}{l}
Y_1 \\ 
Y_2
\end{array}
\right) \right\rangle \right\rangle =\left\langle X_1,Y_1\right\rangle
-\left\langle X_2,Y_2\right\rangle .\quad  \label{f411}
\end{equation}
We are interested in the class of Poisson brackets ${\frak J}$ of the form

\begin{equation}
\left\{ \varphi ,\psi \right\} =\left\langle \left\langle R\left( 
\begin{array}{c}
\nabla \varphi \\ 
\nabla ^{\prime }\varphi
\end{array}
\right) ,\left( 
\begin{array}{c}
\nabla \psi \\ 
\nabla ^{\prime }\psi
\end{array}
\right) \right\rangle \right\rangle  \label{f412}
\end{equation}
where $R\in End{\frak J}\tsize\bigoplus {\frak J}$ , \ $R=$ $-R^{*}.$

\begin{lemma}
\label{l43} The obstruction term in the Jacobi identity for the bracket (\ref
{f412}) is given by 
\begin{equation}
\Delta \equiv \left\{ \left\{ \varphi ,\psi \right\} ,\chi \right\} +{\rm %
c.p.}=\left\langle \left\langle \left[ RD\varphi ,RD\psi \right] ,D\chi
\right\rangle \right\rangle +{\rm c.p.}  \label{f413p}
\end{equation}
\end{lemma}

\begin{theorem}
\label{t42}Let $R\in End\left( {\frak J}\tsize\bigoplus {\frak J}\right) $ be a
skew-symmetric classical r-matrix satisfying the modified classical
Yang-Baxter equation (mCYBE) in ${\frak J}\tsize\bigoplus {\frak J}:$%
\begin{equation*}
\left[ R{\Bbb X},R{\Bbb Y}\right] -R\left( \left[ R{\Bbb X},{\Bbb Y}\right]
+\left[ {\Bbb X},R{\Bbb Y}\right] \right) =-\frac 14\left[ {\Bbb X},{\Bbb Y}%
\right] ,\quad \forall {\Bbb X},{\Bbb Y}\in {\frak J}\tsize\bigoplus {\frak J.}
\end{equation*}
Then the bracket (\ref{f413p}) satisfies the Jacobi identity.
\end{theorem}

{\em Proof.}

\begin{lemma}
\label{l44}If $R$ satisfies mCYBE, then 
\begin{equation}
\left\langle \left\langle \left[ R{\Bbb X},R{\Bbb Y}\right] ,{\Bbb Z}%
\right\rangle \right\rangle +{\rm c.p.}=-\frac 14\left\langle \left\langle
\left[ {\Bbb X},{\Bbb Y}\right] ,{\Bbb Z}\right\rangle \right\rangle .
\label{f416p}
\end{equation}
\end{lemma}

Let $\varphi ,\psi ,\chi $ be some linear functionals on ${\frak J}$. Put 
\begin{equation}
{\Bbb X}=D\varphi ,\quad {\Bbb Y}=D\psi ,\quad {\Bbb Z}=D\chi .  \label{f413}
\end{equation}
By lemmas \ref{l43}, \ref{l44} we have 
\begin{equation*}
\begin{array}{lll}
-4\Delta & = & \left\langle \left\langle \left[ {\Bbb X},{\Bbb Y}\right] ,%
{\Bbb Z}\right\rangle \right\rangle =\left\langle \left[ \nabla \varphi
,\nabla \psi \right] ,\nabla \chi \right\rangle -\left\langle \left[ \nabla
^{\prime }\varphi ,\nabla ^{\prime }\psi \right] ,\nabla ^{\prime }\chi
\right\rangle \\ 
& = & \left\langle \left[ \nabla \varphi ,\nabla \psi \right] ,\nabla \chi
\right\rangle -\limfunc{Tr}d\varphi Ld\psi Ld\chi L+\limfunc{Tr}d\psi
Ld\varphi Ld\chi L \\ 
& = & \left\langle \left[ \nabla \varphi ,\nabla \psi \right] ,\nabla \chi
\right\rangle -\limfunc{Tr}Ld\varphi Ld\psi Ld\chi +\limfunc{Tr}Ld\psi
Ld\varphi Ld\chi =0
\end{array}
\end{equation*}
as desired.

\begin{theorem}
\label{t41}Let ${\frak J}$ be a Lie algebra which is (as a linear space) the
direct sum of the three subalgebras 
\begin{equation}
{\frak J}={\frak J}_{+}\dotplus {\frak J}_0\dotplus {\frak J}_{-},
\label{f41}
\end{equation}
where ${\frak J}_0$ is abelian and 
\begin{equation}
\left[ {\frak J}_{\pm },{\frak J}_0\right] \subset {\frak J}_{\pm }.
\label{f42}
\end{equation}
Let $P_{\pm },$ $P_0$ be the projection operators onto ${\frak J}_{\pm },%
{\frak J}_0$ parallel to the complement and $a,b,c,d$ arbitrary linear
operators in ${\frak J}_0$. Put $r=\frac 12\left( P_{+}-P_{-}\right) .$ Then
the r-matrix 
\begin{equation}
R=\left( 
\begin{array}{cc}
r+aP_0 & bP_0 \\ 
cP_0 & r+dP_0
\end{array}
\right) \in End\left( {\frak J}\tsize\bigoplus {\frak J}\right)  \label{f43}
\end{equation}
satisfies mCYBE.
\end{theorem}

Let us now turn back to the algebra $\Psi {\bf D}_q$ . From theorems \ref
{t42}, \ref{t41} it follows immediately that the q-deformed Gelfand-Dickey
bracket (\ref{f134}) satisfies the Jacobi identity.

\begin{proposition}
\label{p45}The bracket (\ref{f343}) on $\hat G_\alpha ,$ $\alpha $ is
generic, satisfies the Jacobi identity.
\end{proposition}

{\em Proof.} Consider the following bracket on $\Psi {\bf D}_q:$%
\begin{eqnarray}
\left\{ \varphi ,\psi \right\} &=&\left\langle \left\langle \left( 
\begin{array}{cc}
P_{+}+\left( \frac 12+a\right) P_0 & \left( \frac 12-a\right) \hat
h^{-\alpha }P_0 \\ 
\hat h^\alpha \left( \frac 12+a\right) P_0 & P_{+}+\left( \frac 12-a\right)
P_0
\end{array}
\right) \left( 
\begin{array}{c}
\nabla \varphi \\ 
\nabla ^{\prime }\varphi
\end{array}
\right) ,\left( 
\begin{array}{c}
\nabla \psi \\ 
\nabla ^{\prime }\psi
\end{array}
\right) \right\rangle \right\rangle , \notag\\ 
&&\label{f418}
\end{eqnarray}
where, as above, 
\begin{equation}
a=\frac 12\frac{1+\hat h^\alpha }{1-\hat h^\alpha }P_0^{\prime }.
\label{f419}
\end{equation}
Theorems \ref{t41}, \ref{t42} imply the Jacobi identity for (\ref{f418}).

\begin{lemma}
\label{l46} $\hat G_0$ considered as an affine subspace in $\Psi {\bf D}_q$
is a Poisson subspace for (\ref{f418}).
\end{lemma}

{\em Proof. }Obviously, $J_{\left( -\right) }$ is a Poisson subspace. Any
element $A\in J_{\left( -\right) }$ has the form 
\begin{equation}
A=\sum_{i=0}^\infty u_iD^{-i},\quad u_i\in {\Bbb C}\left( \left(
z^{-1}\right) \right) .  \label{f420}
\end{equation}
It is sufficient to prove that functionals of the form 
\begin{equation}
\varphi _f=\limfunc{Tr}Af,\quad f\in {\Bbb C}\left( \left( z^{-1}\right)
\right) ,  \label{f421}
\end{equation}
are central on $\hat G_0,$ i.e., their Poisson brackets with any other
functional vanish identically on $\hat G_0.$ Clearly, $d\varphi _f=f\in J_0,$
hence $\nabla \varphi _f,\nabla ^{\prime }\varphi _f\in J_{\left( -\right)
}. $ So only the $J_0$-components of $\nabla \varphi _f,\nabla ^{\prime
}\varphi _f$ contribute to the bracket (\ref{f418}). Taking into account
that $u_0=1$ on $\hat G_0,$ we get 
\begin{equation}
\left( \nabla \varphi _f\right) _0=\left( \nabla ^{\prime }\varphi _f\right)
_0=f.  \label{f422}
\end{equation}
Then 
\begin{equation}
\left\{ \varphi _f,\psi \right\} =\left\langle \left[ \left( a+\tfrac
12\right) +\left( \tfrac 12-a\right) \hat h^{-\alpha }\right] f,\left(
\nabla \psi \right) _0-\hat h^\alpha \left( \nabla ^{\prime }\psi \right)
_0\right\rangle .  \label{f423}
\end{equation}
But $\left( a+\tfrac 12\right) +\left( \tfrac 12-a\right) \hat h^{-\alpha
}=P_{00},$ hence 
\begin{equation*}
\left\{ \varphi _f,\psi \right\} =\limfunc{Tr}f\cdot \left( \limfunc{Tr}%
\nabla \psi -\limfunc{Tr}\hat h^\alpha \left( \nabla ^{\prime }\psi \right)
_0\right) =0
\end{equation*}
due to the invariance of the inner product.

Consider the map $i:\hat G_\alpha \rightarrow \hat G_0$ defined by $i\left(
L\right) =\bar L.$ Obviously, $i$ is a bijection. It is easy to verify that
the pullback of the Poisson bracket (\ref{f418}) coincides with the bracket (%
\ref{f343}), so the latter satisfies the Jacobi identity.

\end{document}